\documentclass[reqno,12pt]{amsart}
\usepackage{a4wide}

\usepackage{amssymb,graphicx}

\newcommand{\CH}[1][n]{\ensuremath{\mathbb{C}H^{#1}}}
\newcommand{\C}{\ensuremath{\mathbb{C}}}
\newcommand{\R}{\ensuremath{\mathbb{R}}}

\newcommand{\sech}{\ensuremath{\mbox{sech}}}

\newcommand{\enabla}{\ensuremath{\bar{\nabla}}}
\newcommand{\eR}{\ensuremath{\bar{R}}}
\newcommand{\tr}{\ensuremath{\mbox{tr}}}

\newtheorem{theorem}{Theorem}[section]
\newtheorem{proposition}[theorem]{Proposition}
\newtheorem{lemma}[theorem]{Lemma}
\newtheorem{corollary}[theorem]{Corollary}

\begin{document}

\title[Real hypersurfaces with constant principal curvatures]
{\Large Real hypersurfaces\\ with constant principal curvatures\\
in complex hyperbolic spaces}

\author{J{\"u}rgen Berndt, Jos{\'e} Carlos D{\'\i}az-Ramos}

\begin{abstract} We present the classification of all real hypersurfaces in
complex hyperbolic space $\CH$, $n \geq 3$, with three distinct constant
principal curvatures.
\end{abstract}


\maketitle

\thispagestyle{empty}

\section{Introduction}

The aim of submanifold geometry is to understand geometric
invariants of submanifolds and to classify submanifolds according
to given geometric data. In Riemannian geometry, the structure of
a submanifold is encoded in the second fundamental form and its
geometry is controlled by the equations of Gau\ss, Codazzi and
Ricci. The situation simplifies for hypersurfaces, as the Ricci
equation is trivial and the second fundamental form can be written
in terms of a self-adjoint tensor field, the shape operator. The
eigenvalues of the shape operator, the so-called principal
curvatures, are the simplest geometric invariants of a
hypersurface. Two basic problems in submanifold geometry are to
understand the geometry of hypersurfaces for which the principal
curvatures are constant, and to classify them. \'{E}lie Cartan
\cite{Ca38} proved that in spaces of constant curvature a
hypersurface has constant principal curvatures if and only if it
is isoparametric. The classification of isoparametric
hypersurfaces has a long history and over the years many
surprising features have been discovered, see \cite{Th00} for a
survey.

Using the Gau\ss-Codazzi equations, \'{Elie} Cartan \cite{Ca38}
also proved that the number $g$ of distinct principal curvatures
of an isoparametric hypersurface in the real hyperbolic space
${\mathbb R}H^n$ is either 1 or 2. This easily leads to a complete
classification: geodesic hyperspheres, horospheres, totally
geodesic hyperplanes and its equidistant hypersurfaces, tubes
around totally geodesic subspaces of dimension $\geq 1$. As a
consequence, all hypersurfaces in real hyperbolic spaces with
constant principal curvatures are open parts of homogeneous
hypersurfaces.

In this paper we deal with the classification problem of real
hypersurfaces with constant principal curvatures in complex
hyperbolic spaces. We briefly describe the current state of the
problem. Obviously, any homogeneous real hypersurface has constant
principal curvatures. The first author and Tamaru \cite{BT04}
derived recently the complete classification of homogeneous real
hypersurfaces in ${\mathbb C}H^n$. The number $g$ of distinct
principal curvatures of all these homogeneous real hypersurfaces
is either 2,3,4 or 5. No examples are known of real hypersurfaces
with constant principal curvatures in ${\mathbb C}H^n$ which are
not an open part of a homogeneous real hypersurface. It is also
not known whether for any real hypersurface with constant
principal curvatures in ${\mathbb C}H^n$ the number $g$ of
distinct principal curvatures must necessarily be 2,3,4 or 5.

From the Codazzi equation one can easily deduce  that $g > 1$ (see
Corollary \ref{noumbilical}). It follows from work by Montiel
\cite{Mo85} that every real hypersurface with two distinct
constant principal curvatures in complex hyperbolic space
${\mathbb C}H^n$, $n \geq 3$, is an open part of a geodesic
hypersphere, of a horosphere, of a tube around a totally geodesic
${\mathbb C}H^{n-1} \subset {\mathbb C}H^n$, or of a tube with
radius $\ln(2+\sqrt{3})$ around a totally geodesic ${\mathbb R}H^n
\subset {\mathbb C}H^n$. For $n=2$ this problem appears to be
still open. In Corollary \ref{twoconstant} we present a proof for
this classification which includes this low-dimensional case as
well. All these real hypersurfaces are homogeneous Hopf
hypersurfaces. If $\xi$ is a (local) unit normal field of a real
hypersurface $M$ in a Hermitian manifold $\bar{M}$, and $J$
denotes the complex structure of $\bar{M}$, then the Hopf vector
field $J\xi$ is tangent to $M$ everywhere. The hypersurface $M$ is
said to be a Hopf hypersurface if the integral curves of $J\xi$
are geodesics in $M$. If $\bar{M}$ is a K\"{a}hler manifold this
is equivalent to the condition that $J\xi$ is a principal
curvature vector of $M$ everywhere.

The first author obtained in \cite{Be89} the classification of all
Hopf hypersurfaces with constant principal curvatures in ${\mathbb
C}H^n$. Any such hypersurface is an open part of a horosphere, of
a tube around a totally geodesic ${\mathbb C}H^k \subset {\mathbb
C}H^n$ for some $k \in \{0,\ldots,n-1\}$, or to a tube around a
totally geodesic ${\mathbb R}H^n \subset {\mathbb C}H^n$. All
these tubes and the horospheres are homogeneous hypersurfaces and
satisfy $g \in \{2,3\}$. But not all homogeneous real
hypersurfaces in ${\mathbb C}H^n$ are necessarily Hopf
hypersurfaces, see \cite{Be98} for the construction of the
following examples.

Let $KAN$ be an Iwasawa decomposition of $SU(1,n)$, the connected
component of the isometry group of ${\mathbb C}H^n$. The solvable
Lie group $AN$ acts simply transitively on ${\mathbb C}H^n$. The
Riemannian metric on ${\mathbb C}H^n$ therefore induces in a
natural way an inner product on the Lie algebra ${\mathfrak a}
\oplus {\mathfrak n}$ of $AN$. The nilpotent Lie group $N$ is
isomorphic to the $(2n-1)$-dimensional Heisenberg group, and the
orbits of the action of $N$ on ${\mathbb C}H^n$ give a foliation
by horospheres. The Lie algebra ${\mathfrak n}$ of $N$ is a
Heisenberg algebra and has a natural orthogonal decomposition
${\mathfrak n} = {\mathfrak z} \oplus {\mathfrak v}$, where
${\mathfrak z}$ is the one-dimensional center of ${\mathfrak n}$.
Let ${\mathfrak w}$ be a linear hyperplane of ${\mathfrak v}$.
Then ${\mathfrak a} \oplus {\mathfrak z} \oplus {\mathfrak w}$ is
a subalgebra of ${\mathfrak a} \oplus {\mathfrak n}$ of
codimension one. The corresponding connected Lie subgroup of $AN$
therefore induces a foliation on ${\mathbb C}H^n$ by homogeneous
hypersurfaces. None of these homogeneous hypersurfaces is a Hopf
hypersurface. Exactly one of the orbits is minimal and has a
simple geometric description. Consider a totally geodesic
${\mathbb R}H^2 \subset {\mathbb C}H^2 \subset {\mathbb C}H^n$ and
pick a horocycle ${\gamma}$ in ${\mathbb R}H^2$. At each point $p
\in {\gamma}$ we attach the totally geodesic complex hyperbolic
hyperplane which is tangent to the orthogonal complement of the
complex span of the tangent line to ${\gamma}$ at $p$. In this way
we obtain a ruled real hypersurface $W^{2n-1}$ in ${\mathbb
C}H^n$. This hypersurface $W^{2n-1}$ is congruent to the unique
minimal orbit in the above foliation. As can be seen from the
construction, the other homogeneous hypersurfaces in the foliation
are geometrically the equidistant hypersurfaces to $W^{2n-1}$.

It was shown in \cite{Be98} that each of these homogeneous
hypersurfaces has three distinct constant principal curvatures.
Saito claims in \cite{Sa99} that every real hypersurface with
three distinct constant principal curvatures in ${\mathbb C}H^n$
is a Hopf hypersurface, and hence the assumption in \cite{Be89} on
the Hopf hypersurface would be redundant. The above examples show
that this is not true.

The construction of $W^{2n-1}$ can be generalized in the following
way. Consider a totally geodesic ${\mathbb R}H^{k+1} \subset
{\mathbb C}H^n$, $1 \leq k \leq n-1$, and fix a horosphere $H$ in
${\mathbb R}H^{k+1}$. At each point $p \in H$ we attach the
totally geodesic ${\mathbb C}H^{n-k}$ which is tangent to the
orthogonal complement of the complex span of the tangent space to
$H$ at $p$. In this way we obtain a $(2n-k)$-dimensional ruled
minimal submanifold $W^{2n-k}$ in ${\mathbb C}H^n$ with totally
real normal bundle of rank $k$. In terms of the above Iwasawa
decomposition, denote by $o \in {\mathbb C}H^n$ the fixed point of
the action of the compact group $K$ on ${\mathbb C}H^n$. Then
$W^{2n-k}$ is holomorphically congruent to the orbit through $o$
of the closed subgroup of $AN$ with Lie algebra ${\mathfrak a}
\oplus {\mathfrak z} \oplus {\mathfrak w}$, where ${\mathfrak w}$
is the orthogonal complement in ${\mathfrak v}$ of a real subspace
of ${\mathfrak v}$. For $k=1$ we just obtain the above ruled real
hypersurface. For $k > 1$ the tubes around $W^{2n-k}$ are
homogeneous hypersurfaces (see \cite{BB01}) and hence have
constant principal curvatures. The number of distinct principal
curvatures is four except for the radius $r=\ln(2+\sqrt{3})$,
where there are just three distinct principal curvatures.

In this paper we obtain the classification of all real
hypersurfaces in ${\mathbb C}H^n$ with three distinct constant
principal curvatures .

\begin{theorem}\label{thClassification}
Let $M$ be a connected real hypersurface in ${\mathbb C}H^n$, $n
\geq 3$, with three distinct constant principal curvatures. Then
$M$ is holomorphically congruent to an open part of one of the
following real hypersurfaces:
\begin{enumerate}
\item[(a)] the tube of radius $r>0$ around the totally geodesic
${\mathbb C}H^k \subset {\mathbb C}H^n$ for some
$k\in\{1,\dots,n-2\}$;

\item[(b)] the tube of radius $r > 0$, $r\neq\ln(2+\sqrt{3})$,
around the totally geodesic $\mathbb{R}H^n \subset {\mathbb
C}H^n$;

\item[(c)] the ruled minimal real hypersurface $W^{2n-1} \subset
{\mathbb C}H^n$, or to one of the equidistant hypersurfaces to
$W^{2n-1}$;

\item[(d)] the tube of radius $r=\ln(2+\sqrt{3})$ around the ruled
minimal submanifold $W^{2n-k} \subset {\mathbb C}H^n$ for some $k
\in \{2,\ldots,n-1\}$.
\end{enumerate}
\end{theorem}

For $n=2$ the problem remains open. The hypersurfaces in (a) and
(b) are Hopf hypersurfaces, the hypersurfaces in (c) and (d) are
not Hopf hypersurfaces, and all of them are homogeneous. For the
proof, we first derive some rigidity results of the ruled minimal
submanifolds $W^{2n-k}$ in terms of certain geometric data. In
view of the known classification of Hopf hypersurfaces with
constant principal curvatures in $\CH$ (see \cite{Be89}), we may
assume that $M$ is not a Hopf hypersurface. Using the
Gau\ss-Codazzi equations and Jacobi field theory we then show that
one of the focal sets or equidistant hypersurfaces of $M$ has
these geometric data.

We briefly describe the contents of this paper. In Section
\ref{scFormulae} we derive from the Gau\ss-Codazzi equations some
basic formulae for real hypersurfaces in $\CH$ with constant
principal curvatures, and settle the cases $g \leq 2$. The above
mentioned rigidity results for the ruled minimal submanifolds are
proved in Section \ref{ruled}. In Section \ref{princurv} we
determine the principal curvatures and some other geometric data
for real hypersurfaces in $\CH$ with three distinct constant
principal curvatures. Using Jacobi field theory we then proof the
classification result in Section \ref{proofmain}.

The second author has been supported by project BFM 2003-02949
(Spain).

\section{Preliminaries}\label{scFormulae}

We denote by $\CH$ the $n$-dimensional complex hyperbolic space
equipped with the Fubini Study metric $\langle \cdot , \cdot
\rangle$ of constant holomorphic sectional curvature $-1$. We
assume $n \geq 2$ and denote by $\bar{\nabla}$ and $\bar{R}$ the
Levi Civita covariant derivative and the Riemannian curvature
tensor of $\CH$, respectively, using the sign convention
$\bar{R}_{XY} = [\bar{\nabla}_X,\bar{\nabla}_Y] -
\bar{\nabla}_{[X,Y]}$. Then
$$
\bar{R}_{XY}Z = -\frac{1}{4}\Big(\langle Y,Z\rangle X-\langle
X,Z\rangle Y +\langle JY,Z\rangle JX-\langle JX,Z\rangle JY
-2\langle JX,Y\rangle JZ\Big),
$$
where $J$ is the complex structure of $\CH$. We also write
$\bar{R}_{XYZW}=\langle \bar{R}_{XY}Z,W\rangle$.

Let $M$ be a connected submanifold of $\CH$. We denote by $\nabla$
and $R$ the Levi Civita covariant derivative and the Riemannian
curvature tensor of $M$, respectively. By $TM$ and $\nu M$ we
denote the tangent bundle and the normal bundle of $M$,
respectively. By $\Gamma(TM)$ and $\Gamma(\nu M)$ we denote the
module of all vector fields tangent and normal to $M$,
respectively. Let $X,Y,Z,W \in \Gamma(TM)$ and $\xi\in \Gamma(\nu
M)$.

The Levi Civita covariant derivatives of $M$ and $\CH$ are related
by the Gau\ss\ formula
$$
\bar{\nabla}_X Y=\nabla_X Y+I\!I(X,Y),
$$
where $I\!I$ is the second fundamental form of $M$. The Weingarten
formula is
$$
\bar{\nabla}_X\xi = -S_\xi X + \nabla^\perp_X\xi,
$$
where $S_\xi$ denotes the shape operator of $M$ with respect to
$\xi$ and $\nabla^\perp$ is the induced covariant derivative on
$\nu M$. The second fundamental form and shape operator are
related by $\langle S_\xi X,Y \rangle = \langle I\!I(X,Y),\xi
\rangle$. If $M$ is a real hypersurface and $\xi$ is a unit normal
vector field on $M$, we often write $S$ instead of $S_\xi$. The
fundamental equations of second order of interest to us are the
Gau\ss\ equation
$$
\bar{R}_{XYZW}= R_{XYZW} - \langle I\!I(Y,Z),I\!I(X,W)\rangle
    +\langle I\!I(X,Z),I\!I(Y,W)\rangle
$$
and the Codazzi equation
$$
\bar{R}_{XYZ\xi} = \langle (\nabla_X^\perp
I\!I)(Y,Z)-(\nabla_Y^\perp I\!I)(X,Z),\xi \rangle,
$$
where the covariant derivative of the second fundamental form is
given by
$$
(\nabla_X^\perp I\!I)(Y,Z)=\enabla_X^\perp I\!I(Y,Z)
    -I\!I(\nabla_X Y,Z)-I\!I(Y,\nabla_X Z).
$$

If $M$ is a connected real hypersurface of $\CH$ and $\xi$ is a
global unit normal vector field on $M$, the equations simplify to
\begin{eqnarray*}
\bar{\nabla}_X Y&=&\nabla_X Y+\langle SX,Y \rangle \xi,\\
\bar{\nabla}_X\xi &=& -SX,\\
\bar{R}_{XYZW}&=& R_{XYZW} - \langle SY,Z\rangle\langle SX,W
\rangle
    +\langle SX,Z\rangle\langle SY,W\rangle,\\
\bar{R}_{XYZ\xi} &=& \langle (\nabla_XS)Y-(\nabla_YS)X,Z \rangle.
\end{eqnarray*}

We assume from now on that $M$ is a connected real hypersurface of
$\CH$ with constant principal curvatures. For each principal
curvature $\lambda$ of $M$ we denote by $T_\lambda$ the
distribution on $M$ formed by the principal curvature spaces of
$\lambda$.  By $\Gamma(T_\lambda)$ we denote the set of all
sections in $T_\lambda$, that is, all vector fields on $M$
satisfying $SX = \lambda X$.

The Codazzi equation readily implies
\begin{lemma}\label{thThreeEigenv} For all $X\in\Gamma(T_{\lambda_i})$,
$Y\in\Gamma(T_{\lambda_j})$ and $Z\in\Gamma(T_{\lambda_k})$ we
have
$$
\bar{R}_{XYZ\xi} = ({\lambda_j}-{\lambda_k})\langle\nabla_X Y,Z\rangle
    -(\lambda_i-{\lambda_k})\langle\nabla_Y X,Z\rangle.
$$
\end{lemma}

Assume that $\lambda_i={\lambda_j}={\lambda_k}$ in the previous
lemma. Then $\bar{R}_{XYZ\xi}=0$ for all $X,Y,Z \in
\Gamma(T_{\lambda_i})$. Choosing $Z=X$ we get $0=
\langle JX,Y\rangle\langle X,J\xi\rangle$ for all $X,Y
\in \Gamma(T_{\lambda_i})$, which implies $ 0=4\langle
X,J\xi\rangle\,\bar{R}_{XYZ\xi} =\langle
JY,Z\rangle\,\langle X,J\xi\rangle^2$ for all $X,Y,Z \in
\Gamma(T_{\lambda_i})$. Thus we have proved the following

\begin{lemma}\label{real}
If the orthogonal projection of $J\xi_p$ onto $T_{\lambda_i}(p)$
is nonzero at $p \in M$, then $T_{\lambda_i}(p)$ is a real
subspace of $T_p\CH$, that is, $JT_{\lambda_i}(p) \subset
T_{\lambda_i}^\perp(p)$, where $T_{\lambda_i}^\perp(p)$ is the
orthogonal complement of $T_{\lambda_i}(p)$ in $T_p\CH$.
\end{lemma}

This immediately implies

\begin{corollary} \label{noumbilical}
The number $g$ of distinct principal curvatures of $M$ satisfies
$g > 1$.
\end{corollary}

\begin{corollary}{\rm (Montiel \cite{Mo85} for $n \geq 3$)}\label{twoconstant}
Let $M$ be a connected real hypersurface in $\CH$, $n \geq 2$,
with two distinct constant principal curvatures. Then $M$ is
holomorphically congruent to an open part of a horosphere in
$\CH$, or of a geodesic hypersphere in $\CH$, or of a tube around
a totally geodesic ${\mathbb C}H^{n-1} \subset \CH$, or of the
tube with radius $r = \ln(2+\sqrt{3})$ around a totally geodesic
${\mathbb R}H^n \subset \CH$.
\end{corollary}

\proof We just need to prove that $M$ is a Hopf hypersurface. The
result then follows from the classification of real Hopf
hypersurfaces in $\CH$ with constant principal curvatures (see
\cite{Be89}). Let $\lambda_1,\lambda_2$ be the two principal
curvatures of $M$, and assume that there exists a point $p \in M$
such that $J\xi_p = \langle J\xi_p , u_1 \rangle u_1 + \langle
J\xi_p , u_2 \rangle u_2$ with some unit vectors $u_i \in
T_{\lambda_i}(p)$ and $0 \neq \langle J\xi_p , u_i \rangle$.
According to Lemma \ref{real} both $T_{\lambda_1}(p)$ and
$T_{\lambda_2}(p)$ are real, which implies $JT_{\lambda_1}(p)
\subset T_{\lambda_2}(p) \oplus \R\xi_p$ and $JT_{\lambda_2}(p)
\subset T_{\lambda_1}(p) \oplus \R\xi_p$. Since $n\geq 2$ we can
assume $\dim T_{\lambda_1}(p) \geq 2$. Then we have
$J(T_{\lambda_1}(p)\ominus\R u_1)\subset T_{\lambda_2}(p)$, which
implies $\dim T_{\lambda_2}(p)\geq \dim T_{\lambda_1}(p)-1$. But
$u_2 \notin J(T_{\lambda_1}(p)\ominus\R u_1)$ because of $\langle
u_2,J\xi_p\rangle\neq 0$, and thus we have $\dim
T_{\lambda_2}(p)\geq \dim T_{\lambda_1}(p)$. The previous equality
implies $\dim T_{\lambda_2}(p) \geq 2$, and an analogous argument
yields $\dim T_{\lambda_1}(p)\geq \dim T_{\lambda_2}(p)$.
Therefore, $\dim T_{\lambda_1}(p)=\dim T_{\lambda_2}(p)$. This
implies that $\dim T_p M= \dim T_{\lambda_1}(p)+\dim
T_{\lambda_2}(p)$ is even, which contradicts $\dim M = 2n-1$. \qed

Putting $\lambda_i = \lambda_k$ in Lemma \ref{thThreeEigenv} and
then interchanging $Y$ and $Z$ yields

\begin{lemma}\label{thTwoEigenv}
For all $X,Y \in \Gamma(T_{\lambda_i})$ and $Z \in
\Gamma(T_{\lambda_j})$ with $\lambda_i \neq {\lambda_j}$ we have
$$
4({\lambda_j}-{\lambda_i})\langle\nabla_XY,Z\rangle =
\langle JY,Z \rangle
\langle X,J\xi\rangle + \langle JX,Y\rangle \langle Z,J\xi\rangle
+ 2\langle JX,Z\rangle \langle Y,J\xi \rangle.
$$
\end{lemma}

\begin{corollary}\label{thNablaXX}
For all $X\in\Gamma(T_{\lambda_i})$ with $\langle X,J\xi\rangle=0$ we have
$\nabla_X X\in\Gamma(T_{\lambda_i})$.
\end{corollary}

The following equation is a consequence of the Gau\ss\ and Codazzi
equations and will be used later to obtain some relations among
the principal curvatures.

\begin{lemma}\label{thEigenGauss}
For all unit vector fields
$X\in\Gamma(T_{\lambda_i})$ and $Y\in\Gamma(T_{\lambda_j})$
with $\lambda_i\neq \lambda_j$ we have
\begin{eqnarray*}
0&=&({\lambda_j}-{\lambda_i})\big( 1-4{\lambda_i}\,{\lambda_j}
        +2\langle JX,Y\rangle^2
        +8\langle\nabla_X Y,\nabla_Y X\rangle
        -4\langle\nabla_X X,\nabla_Y Y\rangle\big)\\
&& + 4\langle JX,Y\rangle \big(X\langle Y,J\xi\rangle
        +Y\langle X,J\xi\rangle\big)\\
        &&+\langle X,J\xi\rangle \big(3Y\langle JX,Y\rangle
        +\langle \nabla_Y X,JY\rangle - 2\langle \nabla_X Y,JY\rangle\big)\\
        &&+\langle Y,J\xi\rangle \big(3X\langle JX,Y\rangle
        -\langle \nabla_X Y,JX\rangle + 2\langle \nabla_Y X,JX\rangle\big).
\end{eqnarray*}
\end{lemma}

\begin{proof} The Gau\ss\ equation implies
$$
4R_{XYYX}=(4{\lambda_i}{\lambda_j}-1)
-3\langle JX,Y\rangle^2.
$$
On the other hand, the definition of $R$ yields
$$
\begin{array}{rcl}
R_{XYYX} &=& \displaystyle\langle\nabla_X\nabla_Y
Y-\nabla_Y\nabla_X Y
        -\nabla_{[X,Y]}Y,X\rangle\\
\noalign{\medskip} &=& \displaystyle X\langle\nabla_Y Y,X\rangle
        -\langle\nabla_X X,\nabla_Y Y\rangle
        -Y\langle\nabla_X Y,X\rangle\\
\noalign{\medskip} &&  \displaystyle+\langle\nabla_X Y,\nabla_Y
X\rangle
        -\langle\nabla_{[X,Y]}Y,X\rangle.
\end{array}
$$
From Lemma \ref{thTwoEigenv} we get
\begin{eqnarray*}
4({\lambda_j}-{\lambda_i})X\langle\nabla_Y Y,X\rangle
&=& 3\langle Y,J\xi\rangle X\langle JX,Y\rangle
    +3\langle JX,Y\rangle X\langle Y,J\xi\rangle,\\
4({\lambda_i}-{\lambda_j})Y\langle\nabla_X Y,X\rangle
&=& 3\langle X,J\xi\rangle Y\langle JX,Y\rangle
    +3\langle JX,Y\rangle Y\langle X,J\xi\rangle.
\end{eqnarray*}
Next, using the Codazzi equation and the algebraic Bianchi identity, we get
\begin{eqnarray*}
&&(\lambda_j-\lambda_i)\langle\nabla_{[X,Y]}Y,X\rangle\\
&=&        \langle\left(\nabla_{[X,Y]}S\right) Y,X\rangle\\
&=& \langle\left(\nabla_Y S\right)[X,Y], X\rangle +
        \eR_{[X,Y] Y X \xi}\\
&=& \langle\left(\nabla_Y S\right)X, \nabla_X Y\rangle
        - \langle\left(\nabla_Y S\right)X, \nabla_Y X\rangle
       + \eR_{[X,Y] Y X \xi}\\
&=& \langle\left(\nabla_Y S\right)X, \nabla_X Y\rangle
        - \langle\left(\nabla_X S\right)Y, \nabla_Y X\rangle
       - \eR_{Y X \nabla_YX \xi}+ \eR_{[X,Y] Y X\xi}\\
&=& ({\lambda_i}-{\lambda_j})\langle\nabla_XY,\nabla_Y X\rangle
       + \eR_{\nabla_XY Y X \xi}+ \eR_{X \nabla_Y X Y \xi}\\
&=& ({\lambda_i}-{\lambda_j})\langle\nabla_XY,\nabla_Y X\rangle
        \\
&&
\displaystyle-\frac{1}{4}\Big(
        (\lambda_i-\lambda_j)\langle JX,Y\rangle^2
        +\langle JX,Y\rangle \big(X\langle Y,J\xi\rangle
        +Y\langle X,J\xi\rangle\big)\\
&& \displaystyle\phantom{xxxx}
        +\langle X,J\xi\rangle \big(
        \langle JY,\nabla_Y X\rangle - 2\langle JY,\nabla_X Y\rangle\big)\\
&& \displaystyle\phantom{xxxx}
        -\langle Y,J\xi\rangle \big(
        \langle JX,\nabla_X Y\rangle - 2\langle JX,\nabla_Y X\rangle\big)\Big)
\end{eqnarray*}
Altogether this implies the lemma. \end{proof}

\section{The ruled minimal submanifolds $W^{2n-k}$}\label{ruled}

In this section we present a characterization of the ruled minimal
submanifolds $W^{2n-k}$, $k \in \{1,\ldots,n-1\}$. Let $KAN$ be an
Iwasawa decomposition of $SU(1,n)$ and $o \in \CH$ the fixed point
of the action of $K$ on $\CH$. Then $AN$ acts simply transitively
on $\CH$ and we can identify $\CH$ with the solvable Lie group
$AN$ equipped with a suitable left-invariant metric. This induces
an inner product on the Lie algebra ${\mathfrak a} \oplus
{\mathfrak n}$ of $AN$. There is a natural decomposition of the
Lie algebra ${\mathfrak n} = {\mathfrak z} \oplus {\mathfrak v}$
of $N$, where ${\mathfrak z}$ is the one-dimensional center of
${\mathfrak n}$ and ${\mathfrak v}$ is the orthogonal complement
of ${\mathfrak z}$ in ${\mathfrak n}$. The K\"{a}hler structure on
$\CH$ induces a complex structure $i$ on the vector space
${\mathfrak v}$, so that ${\mathfrak v}$ becomes isomorphic to the
complex vector space ${\mathbb C}^{n-1}$. Let ${\mathfrak w}$ be a
linear subspace of ${\mathfrak v}$ such that the orthogonal
complement ${\mathfrak w}^\perp$ of ${\mathfrak w}$ in ${\mathfrak
v}$ is a real subspace of dimension $k$. Then ${\mathfrak s} =
{\mathfrak a} \oplus {\mathfrak z} \oplus {\mathfrak w}$ is a
subalgebra of ${\mathfrak a} \oplus {\mathfrak n}$, and the orbit
through $o$ of the closed subgroup $S$ of $AN$ with Lie algebra
${\mathfrak s}$ is holomorphically congruent to the ruled minimal
submanifold $W^{2n-k}$.

Let ${\mathfrak w}_{\mathbb C}$ be the maximal complex subspace of
${\mathfrak w}$, that is, the orthogonal complement in ${\mathfrak
w}$ of $i{\mathfrak w}^\perp$. Then we have an orthogonal
decomposition ${\mathfrak w} = {\mathfrak w}_{\mathbb C} \oplus
i{\mathfrak w}^\perp$. The subspace ${\mathfrak a} \oplus
{\mathfrak z} \oplus {\mathfrak w}_{\mathbb C}$ is a subalgebra of
${\mathfrak a} \oplus {\mathfrak n}$, and the corresponding Lie
subgroup of $AN$ induces a foliation of $W^{2n-k}$ by totally
geodesic ${\mathbb C}H^{n-k} \subset \CH$. The subspace
${\mathfrak a} \oplus i{\mathfrak w}^\perp$ is a subalgebra of
${\mathfrak a} \oplus {\mathfrak n}$, and the corresponding Lie
subgroup of $AN$ induces a foliation of $W^{2n-k}$ by totally
geodesic ${\mathbb R}H^{k+1} \subset \CH$. Moreover, the subspace
$i{\mathfrak w}^\perp$ is a subalgebra of ${\mathfrak a} \oplus
{\mathfrak n}$, and the corresponding Lie subgroup of $AN$ induces
a foliation of $W^{2n-k}$ by Euclidean spaces ${\mathbb R}^k$
which are embedded in the real hyperbolic spaces ${\mathbb
R}H^{k+1}$ as horospheres in the usual way as totally umbilical
submanifolds with parallel mean curvature vector.

The procedure for the computation of the Levi Civita connection
for a left-invariant Riemannian metric on a Lie group is
well-known and allows us to calculate the second fundamental form
$I\!I$ of $W^{2n-k}$ in an elementary way via the Gau\ss\ formula.
Using for instance the expression for the Levi Civita connection
of $AN$ given in \cite{BTV95}, p.\ 84, shows that $I\!I$ is
determined by
$$
\forall\  \xi \in {\mathfrak w}^\perp : 2I\!I(Z,i\xi) = \xi\ ,
$$
where $Z \in {\mathfrak z}$ is a unit vector with a suitable
orientation, and $Z,\xi,i\xi$ are viewed as left-invariant vector
fields on $AN$. In other words, let $\xi \in {\mathfrak w}^\perp$
be a unit normal vector field of $W^{2n-k}$. Then the principal
curvatures of $W^{2n-k}$ with respect to $\xi$ are
$0,\frac{1}{2},-\frac{1}{2}$ with multiplicities $2n-k-2,1,1$,
respectively, and the principal curvature spaces with respect to
$\pm\frac{1}{2}$ are spanned by $Z \pm i\xi$. This clearly shows
that $W^{2n-k}$ is a minimal submanifold of $\CH$.

We will now show that this second fundamental form characterizes
$W^{2n-k}$ among all $(2n-k)$-dimensional submanifolds of $\CH$
with totally real normal bundle.

\begin{theorem}\label{ruledk}
Let $M$ be a $(2n-k)$-dimensional connected submanifold in $\CH$,
$n \geq 3$, with totally real normal bundle $\nu M \subset T\CH$.
Assume that there exists a unit vector field $Z$ tangent to the
maximal holomorphic subbundle of $TM \subset T\CH$ such that the
second fundamental form $I\!I$ of $M$ is given by the trivial
bilinear extension of $2I\!I(Z,J\xi) = \xi$ for all $\xi \in
\Gamma(\nu M)$. Then $M$ is holomorphically congruent to an open
part of the ruled minimal submanifold $W^{2n-k}$.
\end{theorem}

\begin{proof} We will show the following:
\begin{itemize}
\item[(i)] The maximal holomorphic subbundle ${\mathfrak D}$ of
$TM$ is integrable and each integral manifold is an open part of a
totally geodesic $\C H^{n-k} \subset \CH$; \item[(ii)] The totally
real subbundle ${\mathbb R}JZ \oplus J(\nu M)$ of TM is integrable
and each integral manifold is an open part of a totally geodesic
$\R H^{k+1} \subset \CH$; \item[(iii)] The totally real subbundle
$J(\nu M)$ of $TM$ is integrable and each integral manifold is an
open part of a horosphere in a totally geodesic $\R H^{k+1}
\subset \CH$.
\end{itemize}
The rigidity of totally geodesic submanifolds of Riemannian
manifolds (see e.g.\ \cite{BCO03}, p.\ 230),
and of horospheres in real hyperbolic spaces
(see e.g.\ \cite{BCO03}, pp.\ 24-26), then
implies the assertion.

Ad (i): For $U,V \in \Gamma({\mathfrak D})$ and $\xi \in
\Gamma(\nu M)$ we have $$\langle \nabla_UV,J\xi\rangle = \langle
\bar{\nabla}_UV,J\xi\rangle = -\langle J\bar{\nabla}_UV,\xi\rangle
= -\langle \bar{\nabla}_UJV,\xi\rangle = -\langle
I\!I(U,JV),\xi\rangle = 0$$ and $$\langle
\bar{\nabla}_UV,\xi\rangle = \langle I\!I(U,V),\xi\rangle = 0.$$
This shows that ${\mathfrak D}$ is an autoparallel subbundle of
$TM$ and each integral manifold is a totally geodesic submanifold
of $\CH$. As ${\mathfrak D}$ is a complex subbundle of complex
rank $n-k$, each of these integral manifolds must be an open part
of a totally geodesic $\C H^{n-k} \subset \CH$.

Ad (ii): Let $X \in \Gamma({\mathfrak D} \ominus {\mathbb R}JZ)$
and $\zeta \in \Gamma(\nu M)$ be a local unit normal vector field
of $M$. Using the Gau\ss\ formula, $\bar\nabla J = 0$, the Codazzi
equation, the assumption on $I\!I$, and the explicit expression
for $\bar{R}$ we get
\begin{eqnarray*}
\langle \nabla_{JZ}JZ,X \rangle
& = & \langle \bar{\nabla}_{JZ}JZ,X \rangle
\ =\ \langle \bar{\nabla}_{JZ}JX,Z \rangle
\ =\ \langle \nabla_{JZ}JX,Z \rangle \\
&=& 2 \langle I\!I(\nabla_{JZ}JX,J\zeta),\zeta) \rangle
\ = \ - 2\langle (\nabla_{JZ}I\!I)(JX,J\zeta),\zeta \rangle \\
&=& - 2\langle (\nabla_{J\zeta}I\!I)(JZ,JX),\zeta \rangle -
2\bar{R}_{JZ\,J\zeta\,JX\,\zeta} \ =\  0.
\end{eqnarray*}
For all $\xi,\eta \in \Gamma(\nu M)$ we get
\begin{eqnarray*}
\langle \nabla_{JZ}J\xi,X \rangle &=& \langle
\bar{\nabla}_{JZ}J\xi,X \rangle \ =\ \langle
\bar{\nabla}_{JZ}JX,\xi \rangle \ =\ \langle I\!I(JZ,JX),\xi
\rangle
\ =\ 0,\\
\langle \nabla_{J\xi}JZ,X \rangle
& = & \langle \bar{\nabla}_{J\xi}JZ,X \rangle
\ =\ \langle \bar{\nabla}_{J\xi}JX,Z \rangle
\ =\ \langle \nabla_{J\xi}JX,Z \rangle \\
&=& 2 \langle I\!I(\nabla_{J\xi}JX,J\zeta),\zeta) \rangle
\ = \ - 2\langle (\nabla_{J\xi}I\!I)(JX,J\zeta),\zeta \rangle \\
&=& - 2\langle (\nabla_{JX}I\!I)(J\xi,J\zeta),\zeta \rangle -
2\bar{R}_{J\xi\,JX\,J\zeta\,\zeta}
\ =\  0,\\
\langle \nabla_{J\xi}J\eta,X \rangle &=& \langle
\bar{\nabla}_{J\xi}J\eta,X \rangle \ =\ \langle
\bar{\nabla}_{J\xi}JX,\eta \rangle \ =\ \langle I\!I(J\xi,JX),\eta
\rangle \ =\ 0.
\end{eqnarray*}
Finally, for all $U,V \in \Gamma({\mathbb R}JZ \oplus J(\nu M))$
we obviously have
$$
\langle \bar{\nabla}_UV,\zeta \rangle = \langle I\!I(U,V),\zeta
\rangle = 0.
$$
Altogether this shows that ${\mathbb R}JZ \oplus J(\nu M)$ is
integrable and each integral manifold is a totally geodesic
submanifold of $\CH$. As ${\mathbb R}JZ \oplus J(\nu M)$ is a
totally real subbundle of rank $k+1$, each of these totally
geodesic submanifolds must be an open part of a totally geodesic
$\R H^{k+1} \subset \CH$.

Ad (iii): For all $\xi,\eta \in \Gamma(\nu M)$ we get
\begin{eqnarray*}
\langle \nabla_{J\xi}J\eta,JZ \rangle
& =& \langle \bar{\nabla}_{J\xi}J\eta,JZ \rangle
\ = \ \langle \bar{\nabla}_{J\xi}\eta,Z \rangle
\ = \ -\langle \bar{\nabla}_{J\xi}Z,\eta \rangle\\
&=& - \langle I\!I(Z,J\xi),\eta \rangle \ =\ -\frac{1}{2} \langle
\xi,\eta \rangle \ =\ -\frac{1}{2} \langle J\xi,J\eta \rangle.
\end{eqnarray*}
It follows that $\langle [J\xi,J\eta],JZ \rangle = 0$ for all
$\xi,\eta \in \Gamma(\nu M)$. Together with (ii) this implies that
$J(\nu M)$ is integrable and the second fundamental form
$\tilde{I\!I}$ of an integral manifold is given by
$$
\tilde{I\!I}(J\xi,J\eta) = - \frac{1}{2}\langle J\xi,J\eta \rangle
JZ.
$$
Thus each integral manifold is a totally umbilical submanifold
with constant mean curvature $1/2$ in a real hyperbolic space $\R
H^{k+1}$ of constant sectional curvature $-1/4$. If $k \geq 2$,
the classification of totally umbilical submanifolds in real
hyperbolic spaces shows that each integral manifold is an open
part of a horosphere in $\R H^{k+1}$. If $k = 1$, we have $
2\bar\nabla_{J\xi}J\xi = 2J\bar\nabla_{J\xi}\xi = -2JSJ\xi = -JZ $
and hence $ 4 \bar\nabla_{J\xi}\bar\nabla_{J\xi}J\xi =  -
2\bar\nabla_{J\xi}JZ = - J\xi $. Thus the integral curves of
$J\xi$ satisfy the differential equation for a horocycle in $\R
H^2$, which implies that the integral manifolds of the
distribution $J(\nu M)$ are open parts of horocycles in $\R
H^2$.\end{proof}

For $k=1$ we have the following improvement:

\begin{theorem}\label{ruled1}
Let $M$ be a connected real hypersurface in $\CH$, $n \geq 3$,
with three distinct principal curvatures $0$, $+1/2$ and $-1/2$
and corresponding multiplicities $2n-3$, $1$ and $1$,
respectively. Then $M$ is holomorphically congruent to an open
part of the ruled real hypersurface $W^{2n-1}$.
\end{theorem}

\begin{proof} Let $p \in M$ and suppose that the orthogonal projection of
$J\xi_p$ onto $T_0(p)$ is nonzero. Then $T_0(p)$ is a real
subspace of $T_p\CH$ by Lemma \ref{real}. Since $\dim T_0(p) =
2n-3$, this is impossible for $n > 3$ and we must have $n=3$.
Since $\xi_p \in T_0^\perp(p)$ it follows that $J\xi_p \in
T_0(p)$. Since orthogonal projection onto subbundles is a
continuous mapping, this must hold on an open neighborhood $U$ of
$p$ in $M$. Therefore, $U$ is a Hopf hypersurface in $\CH[3]$ with
three distinct constant principal curvatures $0$, $+1/2$ and
$-1/2$. According to the classification in \cite{Be89} of Hopf
hypersurfaces with constant principal curvatures in $\CH$ such a
hypersurface does not exist. We conclude that the orthogonal
projection of the Hopf vector field $J\xi$ onto $T_0$ is zero
everywhere.

Now define $M^+$ as the set of all points $p \in M$ at which the
orthogonal projections of $J\xi_p$ onto $T_{1/2}(p)$ and
$T_{-1/2}(p)$ are both nonzero. Clearly, $M^+$ is an open subset
of $M$. Using again the classification in \cite{Be89} of Hopf
hypersurfaces with constant principal curvatures in $\CH$, we see
that $M^+$ is nonempty.

Let $X$ and $Y$ be local unit vector fields on $M$ with $X \in
\Gamma(T_{1/2})$ and $Y \in \Gamma(T_{-1/2})$. Then we can write
$J\xi = aX+bY$ with $a,b \in \R$ such that $a^2+b^2=1$. We may
assume that $X$ and $Y$ are chosen such that $a,b \geq 0$. As we
have seen above, $T_0(p)$ cannot be a real subspace at any point
$p \in M$. Thus there exist vector fields $U,V \in \Gamma(T_0)$
with $\langle JU,V\rangle\neq 0$. Since $\bar{\nabla}J = 0$ we
have $\bar{\nabla}_U J\xi = J\bar{\nabla}_U\xi = -JSU = 0$, and
thus Lemma \ref{thTwoEigenv} implies
$$
0 = U\langle V,J\xi\rangle=\langle\nabla_U V,J\xi\rangle =
a\langle\nabla_U V,X\rangle + b\langle \nabla_UV,Y\rangle =
\frac{1}{2}(a^2-b^2)\langle JU,V \rangle\ .$$ As $\langle
JU,V\rangle\neq 0$ this gives $a^2 = b^2$ and hence $a = b =
1/\sqrt{2}$. This shows that $M^+$ is a closed subset of $M$. As
$M^+$ is open and nonempty, we see that $M^+ = M$. In particular,
the length of the orthogonal projections of the Hopf vector field
$J\xi$ onto $T_{1/2}$ and $T_{-1/2}$ is constant and equal to
$1/\sqrt{2}$. We now define $Z = a(X-Y)$. Then the second
fundamental form of $M$ is of the form as in Theorem \ref{ruledk},
and the result now follows from that theorem.\end{proof}

\section{Principal curvatures}\label{princurv}

Let $M$ be an orientable connected real hypersurface of $\CH$ and
$\xi$ a global unit normal vector field on $M$. We assume that $M$
has three distinct constant principal curvatures $\lambda_1$,
$\lambda_2$, $\lambda_3$ and denote by $m_i$ the multiplicity of
$\lambda_i$. If $M$ is a Hopf hypersurface, it was shown in
\cite{Be89} that $M$ is an open part of a tube around a totally
geodesic $\CH[k]\subset\CH$ for some $k\in\{1,\dots,n-2\}$, or of
a tube with radius $r \neq \ln(2+\sqrt{3})$ around a totally
geodesic $\R H^n\subset\CH$. We can therefore assume that $M$ is
not a Hopf hypersurface. Then there exists an open subset of $M$
on which at least two of the three orthogonal projections of the
Hopf vector field $J\xi$ onto the principal curvature
distributions $T_{\lambda_i}$ are nontrivial.

In the first part of this section we will prove that there cannot be
three nontrivial projections. We then derive some equations
relating the principal curvatures and obtain some geometrical
information about the principal curvature distributions.

\begin{lemma}\label{thJui}
Assume that there exists a point $p\in M$ such that $J\xi_p=\sum
b_i u_i$ with some unit vectors $u_i\in T_{\lambda_i}(p)$ and
$0\neq b_i\in\R$, $i=1,2,3$. Then $\R \xi_p\oplus \R u_1\oplus\R
u_2\oplus\R u_3$ is a complex subspace of $T_p\CH$ and, by a
suitable orientation of $u_1,u_2,u_3$, we have $b_i=\langle
Ju_j,u_k\rangle$ for all cyclic permutations $(i,j,k)$ of
$(1,2,3)$.
\end{lemma}

\begin{proof} According to Lemma \ref{real}, each
of the three principal curvature spaces
$T_{\lambda_i}(p)$ is a real subspace of $T_p{\mathbb C}H^n$.
Thus we can write
\begin{equation}\label{eqJuiSum}
J u_i=\sum_{j=1}^3\langle J u_i,u_j\rangle u_j+\sum_{j=1}^3
w_{ij}-b_i\xi_p,
\end{equation}
with some vectors $w_{ij}\in T_{\lambda_j}(p)\ominus\R u_j$,
$w_{ii}=0$. Then we have
\begin{equation}\label{eqSystem}
0=\langle u_i,\xi_p\rangle=\langle J u_i,J\xi_p\rangle
=\sum_{j=1}^3 b_j\langle J u_i,u_j\rangle,
\end{equation}
and hence
\[
-\xi_p=J^2\xi_p=J(J\xi_p)=\sum_{i=1}^3 b_i J u_i
=\sum_{j=1}^3\left(\sum_{i=1}^3 b_i w_{ij}\right)-\xi_p.
\]
This implies $\sum_{i=1}^3 b_i w_{ij}=0$ for all $j\in\{1,2,3\}$.
Thus for each $j\in\{1,2,3\}$ the two vectors $w_{ij}$ with $i\neq
j$ are either both zero, or both nonzero and collinear. From
(\ref{eqJuiSum}) and $J\xi_p=\sum b_i u_i$, we therefore see that
$\R u_1\oplus\R u_2\oplus\R u_3\oplus\R w_{12}\oplus\R
w_{23}\oplus\R w_{31}\oplus\R \xi_p$ is a complex subspace of
$T_p\CH$. As the real dimension of a complex vector space is even,
at least one of the three vectors $w_{12}$, $w_{23}$, $w_{31}$
must be zero, say $w_{23}=0$, which implies also $w_{13}=0$.
Moreover, for dimension reasons, the vectors $w_{12}$, $w_{31}$
are either both zero or both nonzero. Then, using
(\ref{eqJuiSum}), we get
\[
\begin{array}{rcccccl}
0&=&\langle u_1,w_{12}\rangle&=&\langle Ju_1,Jw_{12}\rangle
    &=&-\langle Ju_1,u_3\rangle\langle w_{32},w_{12}\rangle,\\
\noalign{\medskip}
0&=&\langle u_2,w_{21}\rangle&=&\langle
Ju_2,Jw_{21}\rangle
    &=&-\langle Ju_2,u_3\rangle\langle w_{31},w_{21}\rangle.
\end{array}
\]
If $w_{12}$, $w_{31}$ are both nonzero, then $w_{32}$, $w_{21}$
are nonzero as well, and we get $\langle Ju_1,u_3\rangle=0=\langle
Ju_2,u_3\rangle$ using the collinearity of $w_{12}$, $w_{32}$ and
$w_{31}$, $w_{21}$. From (\ref{eqSystem}) we then get $\langle
Ju_1,u_2\rangle=0$ as well. This implies $Ju_1=w_{12}-b_1\xi_p$
and hence $b_1 J\xi_p=u_1+Jw_{12}$. As
$T_{\lambda_2}(p)$ is a real subspace of $T_p\CH$,
the previous equation shows that $J\xi_p\in
T_{\lambda_1}(p)\oplus T_{\lambda_3}(p)$, which contradicts the
assumption on $J\xi_p$. Hence $w_{12}$, $w_{31}$ are both zero.
Altogether this shows that $\R u_1\oplus\R u_2\oplus\R u_3\oplus\R
\xi_p$ is a complex subspace of $T_p\CH$.

Finally, solving the system of equations (\ref{eqSystem}),
we see that the vector
$(b_1,b_2,b_3)$ is in the real span of $(\langle
Ju_2,u_3\rangle,\langle Ju_3,u_1\rangle,\langle Ju_1,u_2\rangle)$.
From $b_1^2+b_2^2+b_3^2=1$ and (\ref{eqJuiSum}) we get
$$
3 = \sum_{i=1}^3\langle u_i,u_i\rangle^2
=\sum_{i=1}^3\langle Ju_i,Ju_i\rangle^2
=2\Bigl(\langle Ju_2,u_3\rangle^2+\langle Ju_3,u_1\rangle^2
        +\langle Ju_1,u_2\rangle^2\Bigr)+1.
$$
Thus $(\langle Ju_2,u_3\rangle,\langle Ju_3,u_1\rangle, \langle
Ju_1,u_2\rangle)$ is a unit vector in $\R^3$, and the lemma now
follows.\end{proof}

\begin{lemma}\label{thRelation3eigenvalues}
Assume that there exists a point $p \in M$ such that the
orthogonal projections of $J\xi_p$ onto $T_{\lambda_i}(p)$,
$i=1,2,3$, are nontrivial. Then we have
\[
(2\lambda_i(\lambda_i-\lambda_j)-1)\langle Jw_i,w_j\rangle=0
\]
for all distinct $i,j\in\{1,2,3\}$, $w_i\in
T_{\lambda_i}(p)\ominus\R u_i$ and $w_j\in
T_{\lambda_j}(p)\ominus\R u_j$.
\end{lemma}

\begin{proof}
By continuity, the orthogonal projections of $J\xi$ onto
$T_{\lambda_i}$, $i=1,2,3$, must be nontrivial on an open
neighborhood of $p$ in $M$. The following calculations hold on
this open neighborhood. It follows from Lemma \ref{thJui} that
there exist unit vector fields $U_\nu \in \Gamma(T_{\lambda_\nu})$
such that $J\xi = \sum b_\nu U_\nu$ with $b_\nu = \langle
JU_\mu,U_\rho \rangle$, where $(\nu,\mu,\rho)$ is a cyclic
permutation of $(1,2,3)$. We note that ${\mathfrak D} = TM \ominus
(\R U_1 \oplus \R U_2 \oplus \R U_3)$ is a $J$-invariant subbundle
of $TM$ by Lemma \ref{thJui}.

In the following we denote by $W_\nu$ and $\widetilde{W}_\nu$
vector fields with values in $T_{\lambda_\nu} \cap {\mathfrak D} =
T_{\lambda_\nu}\ominus\R U_\nu$, $\nu =1,2,3$. Using $\bar\nabla J
= 0$ and the Weingarten formula we get $\bar\nabla_{W_i}J\xi =
J\bar\nabla_{W_i}\xi = -JSW_i = -\lambda_i JW_i$. Since
$T_{\lambda_i}$ is real, and using Lemma \ref{thTwoEigenv}, we get
\[
0=W_i\langle\widetilde{W}_i,J\xi\rangle
    =\sum_{\nu=1}^3 b_\nu\langle\enabla_{W_i}\widetilde{W}_i,U_\nu\rangle
    +\langle\widetilde{W}_i,\enabla_{W_i}J\xi\rangle
    =b_i\langle\nabla_{W_i}\widetilde{W}_i,U_i\rangle.
\]
Hence, $\langle\nabla_{W_i}\widetilde{W}_i,U_i\rangle=0$. As
$T_{\lambda_j}$ is real and ${\mathfrak D}$ is complex, we can
write $JW_j=\widetilde{W}_i+\widetilde{W}_k$ with $k \neq i,j$.
Then, using $\langle\nabla_{W_i}\widetilde{W}_i,U_i\rangle=0$,
Lemma \ref{thTwoEigenv}, and the fact that $T_{\lambda_i}$ is
real, we get
\begin{equation}\label{eqNablaWiJUiWj}
\langle\enabla_{W_i}JU_i,W_j\rangle
    =-\langle\nabla_{W_i}U_i,JW_j\rangle
    =\langle\nabla_{W_i}\widetilde{W}_i,U_i\rangle
    -\langle\nabla_{W_i}U_i,\widetilde{W}_k\rangle=0.
\end{equation}
Next, Lemma \ref{thTwoEigenv} implies
\begin{eqnarray*}
0   &=& W_i\langle W_j,J\xi\rangle
        =\sum_\nu b_\nu\langle\bar\nabla_{W_i}W_j,U_\nu\rangle
        +\langle W_j,\enabla_{W_i}J\xi\rangle\\
    &=& \frac{b_i^2}{2(\lambda_i-\lambda_j)}\langle JW_i,W_j\rangle
        +\sum_{\nu\neq i}b_\nu\langle\nabla_{W_i}W_j,U_\nu\rangle
        -\lambda_i\langle JW_i,W_j\rangle\\
    &=& \frac{b_i^2-2\lambda_i(\lambda_i-\lambda_j)}
            {2(\lambda_i-\lambda_j)}\langle JW_i,W_j\rangle
            +b_j\langle\nabla_{W_i}W_j,U_j\rangle
            +b_k\langle\nabla_{W_i}W_j,U_k\rangle.
\end{eqnarray*}
On the other hand, replacing $JU_i$ by $\sum_\nu \langle JU_i ,
U_\nu \rangle U_\nu - b_i\xi$ and using (3) we get $0 = W_i\langle
W_j,JU_i\rangle =
    \sum_\nu\langle JU_i,U_\nu\rangle
            \langle\nabla_{W_i}W_j,U_\nu\rangle $
and hence
$$
0 =  b_k\langle\nabla_{W_i}W_j,U_j\rangle
        -b_j\langle\nabla_{W_i}W_j,U_k\rangle.
$$
The last two equations provide a system of linear equations with
unknowns $\langle\nabla_{W_i}W_j,U_j\rangle$ and
$\langle\nabla_{W_i}W_j,U_k\rangle$. This linear system has a
unique solution which is given by
\begin{equation}\label{thNablaWiWjUk}
\langle\nabla_{W_i}W_j,U_\nu\rangle
    =\frac{b_\nu(2\lambda_i(\lambda_i-\lambda_j)-b_i^2)}
    {2(\lambda_i-\lambda_j)(1-b_i^2)}\langle JW_i,W_j\rangle\ \
    (\nu\neq i).
\end{equation}
As $T_{\lambda_i}$ is real, we have $\langle
JW_i,\widetilde{W}_k\rangle=\langle JW_i,JW_j\rangle=0$, and using
Lemma \ref{thTwoEigenv} and equation (\ref{thNablaWiWjUk}) (with
$j$ and $k$ interchanged) we get
\[
\langle \bar\nabla_{W_i}W_j,JU_k\rangle
    =\langle\nabla_{W_i}U_k,JW_j\rangle
    =-\langle\nabla_{W_i}\widetilde{W}_i,U_k\rangle
    -\langle\nabla_{W_i}\widetilde{W}_k,U_k\rangle=0.
\]
Replacing $JU_k$ by $\sum_\nu \langle JU_k , U_\nu \rangle U_\nu -
b_k\xi$, this implies
$$
0  = \langle\bar\nabla_{W_i}W_j,JU_k\rangle
        =\sum_\nu\langle\nabla_{W_i}W_j,U_\nu\rangle
        \langle JU_k,U_\nu\rangle,
$$
from which we easily get
\[
\langle\nabla_{W_i}W_j,U_j\rangle=\frac{b_j}{2(\lambda_i-\lambda_j)}
\langle JW_i,W_j\rangle
\]
by using Lemma \ref{thTwoEigenv} once again. By comparison of this
equation with equation (\ref{thNablaWiWjUk}) for $\nu=j$ we
eventually get the result.
\end{proof}

\begin{proposition}\label{no3proj}
If $n\geq 3$, then there exists no point $p \in M$ such that the
orthogonal projections of $J\xi_p$ onto $T_{\lambda_i}(p)$,
$i=1,2,3$, are nontrivial.
\end{proposition}

\begin{proof} As $n \geq 3$, the complex vector space ${\mathfrak D}_p =
\bigoplus_i(T_{\lambda_i}(p)\ominus\R u_i)$ has dimension $\geq
1$. Since each $T_{\lambda_i}(p)\ominus\R u_i$ is real, there
exist $i\neq j$ such that  $\langle Jw_i,w_j\rangle\neq 0$ for
some $w_i\in T_{\lambda_i}(p)\ominus\R u_i$, $w_j\in
T_{\lambda_j}(p)\ominus\R u_j$. From Lemma
\ref{thRelation3eigenvalues} we therefore get
$2\lambda_i(\lambda_i-\lambda_j)-1=2\lambda_j(\lambda_j-\lambda_i)-1=0$,
and thus $\lambda_i^2=\lambda_j^2=1/4$. This argument shows that
$T_{\lambda_k}(p)\ominus\R u_k$ must be trivial, that is, the
third eigenvalue $\lambda_k$ has multiplicity one. Since the
eigenspaces are real it also implies that
$J(T_{\lambda_i}(p)\ominus\R u_i)=T_{\lambda_j}(p)\ominus\R u_j$.

Let $W_i\in\Gamma(T_{\lambda_i}\ominus\R U_i)$ be a unit vector
field which is defined in an open neighborhood of $p$ in $M$, and
define $W_j=JW_i\in\Gamma(T_{\lambda_j}\ominus\R U_j)$, where
$U_\nu$ is as in the previous proof. Applying Lemma
\ref{thEigenGauss} to $W_i$ and $W_j$, and using Corollary
\ref{thNablaXX}, we obtain
\[
3-4\lambda_i\lambda_j
+8\langle\nabla_{W_i}W_j,\nabla_{W_j}W_i\rangle=0.
\]
We have $\nabla_{W_i}W_j\in\Gamma(T_{\lambda_j}\oplus\R
U_i\oplus\R U_k)$ and
$\nabla_{W_j}W_i\in\Gamma(T_{\lambda_i}\oplus \R U_j \oplus\R
U_k)$ by Lemma \ref{thTwoEigenv}, and therefore
$\langle\nabla_{W_i}W_j,\nabla_{W_j}W_i\rangle
=\sum_\nu\langle\nabla_{W_i}W_j,U_\nu\rangle
\langle\nabla_{W_j}W_i,U_\nu\rangle$. The latter sum can be
calculated easily by using Lemma \ref{thTwoEigenv} and equation
(\ref{thNablaWiWjUk}). Using the fact that
$4\lambda_i^2=4\lambda_j^2=1$ this gives
$4\langle\nabla_{W_i}W_j,\nabla_{W_j}W_i\rangle=1$. Inserting this
into the above equation yields $5-4\lambda_i\lambda_j=0$. From
$4\lambda_i^2=4\lambda_j^2=1$ and $\lambda_i \neq \lambda_j$ we
know that $4\lambda_i\lambda_j = -1$, which gives a contradiction.
Therefore there exists no point $p \in M$ such that the orthogonal
projections of $J\xi_p$ onto $T_{\lambda_i}(p)$ are nontrivial.
\end{proof}

\begin{lemma}\label{thJuiJa}
Assume that there exists a point $p\in M$ such that $J\xi_p=b_1
u_1+b_2 u_2$ with some unit vectors $u_i\in T_{\lambda_i}(p)$ and
$0\neq b_i\in\R$, $i=1,2$. Then there exists a unit vector $a\in
T_{\lambda_3}(p)$ such that $\R\xi_p\oplus\R u_1\oplus\R
u_2\oplus\R a$ is a complex subspace of $T_p\CH$ and, by a
suitable orientation of $a$, we have $b_1=\langle Ju_2,a\rangle$,
$b_2=-\langle Ju_1,a\rangle$, $\langle Ju_1,u_2\rangle=0$ and
$Ju_i=(-1)^i b_j a-b_i \xi$ with distinct $i,j\in\{1,2\}$.
\end{lemma}

\begin{proof} The eigenspaces $T_{\lambda_1}(p)$ and $T_{\lambda_2}(p)$
are real subspaces of $T_p\C H^n$ by Lemma \ref{real}.  Therefore
we can write
\begin{eqnarray*}
Ju_1    &=& \langle Ju_1,u_2\rangle u_2+w_{12}+w_{13}-b_1\xi,\\
Ju_2    &=& \langle Ju_2,u_1\rangle u_1+w_{21}+w_{23}-b_2\xi,
\end{eqnarray*}
with $w_{21}\in T_{\lambda_1}(p)\ominus\R u_1$, $w_{12}\in
T_{\lambda_2}(p)\ominus\R u_2$, and $w_{13},w_{23} \in
T_{\lambda_3}(p)$. Hence,
\[
\begin{array}{@{}r@{\,}c@{\,}l@{}}
-\xi_p  &=& J^2\xi_p=J(J\xi_p)=b_1 Ju_1+b_2 Ju_2\\
\noalign{\medskip}
        &=& b_2\langle Ju_2,u_1\rangle u_1
            +b_1\langle Ju_1,u_2\rangle u_2+b_2 w_{21}+b_1
            w_{12}+(b_1 w_{13}+b_2 w_{23})-\xi_p.
\end{array}
\]
This shows that $\langle Ju_1,u_2\rangle=0$, $w_{12}=w_{21}=0$ and
$b_1 w_{13}+b_2 w_{23}=0$. As $b_1,b_2\neq 0$, the vectors
$w_{13}$, $w_{23}$ are either both zero or both nonzero. If
$w_{13}=w_{23}=0$, then $Ju_1=-b_1\xi_p$ and $Ju_2=-b_2\xi_p$,
which is impossible. Hence $w_{13}$, $w_{23}$ are both nonzero and
collinear. Let $a$ be a unit vector in $\R w_{13}=\R w_{23}\subset
T_{\lambda_3}(p)$. Since $Jw_{13}=b_1 J\xi_p-u_1\in \R u_1\oplus\R
u_2$ we get $Ja\in\R u_1\oplus\R u_2$, which shows that
$\R\xi_p\oplus\R u_1\oplus\R u_2\oplus\R a$ is a complex subspace
of $T_p\CH$. The two vectors $Ja,J\xi_p\in\R u_1\oplus\R u_2$ are
orthonormal and $J\xi_p=b_1 u_1+b_2 u_2$. Therefore, by a suitable
orientation of $a$, we can write $Ja=b_2 u_1-b_1 u_2$. As $Ja =
\langle Ja,u_1 \rangle u_1 + \langle Ja , u_2 \rangle u_2$, the
result now follows. \end{proof}

In view of Proposition \ref{no3proj} we can assume from now on
that there exists an open subset of $M$ on which the orthogonal
projection of $J\xi$ onto $T_{\lambda_3}$ is trivial. The
following calculations are done on this open subset. It follows
from Lemma \ref{thJuiJa} that there exist unit vector fields $U_1
\in \Gamma(T_{\lambda_1})$, $U_2 \in \Gamma(T_{\lambda_2})$ and $A
\in \Gamma(T_{\lambda_3})$ such that
$$
\begin{array}{c}
J\xi = b_1U_1+b_2U_2\ , \ JU_i=(-1)^i b_j A-b_i \xi\  (i\neq j)\ ,\\
\noalign{\medskip}
b_1=\langle JU_2,A\rangle\ ,\ b_2=\langle JA,U_1\rangle\ , \
\langle JU_1,U_2\rangle=0.
\end{array}
$$
Below we will use these
relations frequently without referring to them explicitly.
Moreover,
$${\mathfrak D} = TM \ominus (\R U_1 \oplus \R U_2 \oplus \R A) =
(T_{\lambda_1} \ominus \R U_1) \oplus (T_{\lambda_2} \ominus \R
U_2)\oplus (T_{\lambda_3} \ominus \R A)$$ is a $J$-invariant
subbundle.

\begin{lemma}
For $i,j \in \{1,2\}$ with $i \neq j$ we have
\begin{eqnarray}
\nabla_{U_i}U_i &=& (-1)^i\frac{3b_1 b_2}
    {4(\lambda_3-\lambda_i)}A,\label{eqNablaUiUi}\\
\nabla_{U_i}U_j &=& (-1)^j\left(\lambda_i
    +\frac{3b_i^2}{4(\lambda_3-\lambda_i)}\right)A,\label{eqNablaUiUj}\\
\nabla_{U_i}A   &=& (-1)^j\frac{3b_1 b_2}{4(\lambda_3-\lambda_i)}U_i
    +(-1)^i\left(\lambda_i
    +\frac{3b_i^2}{4(\lambda_3-\lambda_i)}\right)
    U_j,\label{eqNablaUiA}\\
\nabla_AU_i     &=& \frac{(-1)^j}{\lambda_i-\lambda_j}\left(
    \frac{b_i^2-2b_j^2}{4}+(\lambda_j-\lambda_3)
    \Bigl(\lambda_i+\frac{3b_i^2}{4(\lambda_3-\lambda_i)}
    \Bigr)\right)U_j,\label{eqNablaAUi}\\
\nabla_AA   &=& 0\label{eqNablaAA}.
\end{eqnarray}
\end{lemma}

\begin{proof} Let $W_i \in \Gamma(T_{\lambda_i} \ominus \R U_i)$, $W_j
\in \Gamma(T_{\lambda_j} \ominus \R U_j)$ and $W_3 \in
\Gamma(T_{\lambda_3} \ominus \R A)$. Since $U_i$ has constant
length, we have $\langle\nabla_{U_i}U_i,U_i\rangle=0$. From Lemma
\ref{thTwoEigenv} we easily get
\[
\langle\nabla_{U_i}U_i,U_j\rangle=\langle\nabla_{U_i}U_i,W_j\rangle=
\langle\nabla_{U_i}U_i,W_3\rangle=0\ ,\
\langle\nabla_{U_i}U_i,A\rangle
=(-1)^i\frac{3b_ib_j}{4(\lambda_3-\lambda_i)}A.
\]
As $T_{\lambda_i}$ is real, we have $\langle
W_i,\enabla_{U_i}J\xi\rangle = \langle
W_i,J\enabla_{U_i}\xi\rangle = -\lambda_i\langle W_i,JU_i\rangle =
0$, and using Lemma \ref{thTwoEigenv} once again we then get
\begin{eqnarray*}
0 &=& U_i\langle W_i,J\xi\rangle
        \ =\ \langle\enabla_{U_i}W_i,J\xi\rangle
        +\langle W_i,\enabla_{U_i}J\xi\rangle\\
    &=& b_i\langle\nabla_{U_i}W_i,U_i\rangle
        +b_j\langle\nabla_{U_i}W_i,U_j\rangle
    \ =\ -b_i\langle\nabla_{U_i}U_i,W_i\rangle.
\end{eqnarray*}
Since $b_i\neq 0$, this implies
$\langle\nabla_{U_i}W_i,U_i\rangle=0$, and equation
(\ref{eqNablaUiUi}) now follows.

Since $U_j$ has constant length, we have
$\langle\nabla_{U_i}U_j,U_j\rangle=0$, from (\ref{eqNablaUiUi}) we
get $\langle\nabla_{U_i}U_j,U_i\rangle=0$, and Lemma
\ref{thTwoEigenv} implies
$\langle\nabla_{U_i}U_j,W_i\rangle=-\langle\nabla_{U_i}W_i,U_j\rangle=0$.
Let $\nu \in \{j,3\}$. Using (\ref{eqNablaUiUi}) and $\langle
W_\nu,\enabla_{U_i}J\xi\rangle = \langle
W_\nu,J\enabla_{U_i}\xi\rangle = -\lambda_i\langle
W_\nu,JU_i\rangle = 0$ we get
\[
0=U_i\langle
W_\nu,J\xi\rangle=\langle\enabla_{U_i}W_\nu,J\xi\rangle+\langle
W_\nu,\enabla_{U_i}J\xi\rangle=b_j\langle\nabla_{U_i}W_\nu,U_j\rangle,
\]
which gives $\langle\nabla_{U_i}U_j,W_\nu\rangle =0$. Finally, $0
= U_i\langle JU_i,U_j\rangle = \langle \bar\nabla_{U_i}JU_i,U_j
\rangle + \langle JU_i, \bar\nabla_{U_i}U_j \rangle =
-\langle\bar\nabla_{U_i}U_i,JU_j \rangle + \langle JU_i,
\bar\nabla_{U_i}U_j \rangle$. Replacing now $JU_i$ and $JU_j$ by
the corresponding expressions in terms of $A$ and $\xi$ we obtain
\[
0=b_j\left(\frac{3b_i^2}{4(\lambda_3-\lambda_i)}
+\lambda_i+(-1)^i\langle\nabla_{U_i}U_j,A\rangle\right).
\]
Altogether this now implies equation (\ref{eqNablaUiUj}).

Since $A$ has constant length, we have
$\langle\nabla_{U_i}A,A\rangle=0$. For $\nu \in \{1,2,3\}$ we get
$0=U_i\langle JU_i,W_\nu\rangle
=\langle\enabla_{U_i}JU_i,W_\nu\rangle +\langle
JU_i,\enabla_{U_i}W_\nu\rangle =
-\langle\enabla_{U_i}U_i,JW_\nu\rangle +\langle
JU_i,\enabla_{U_i}W_\nu\rangle$. The first term vanishes because
of equation (\ref{eqNablaUiUi}). For the second term we replace
$JU_i$ by $(-1)^i b_j A-b_i \xi$, which leads to $0 =
\langle\nabla_{U_i}A,W_\nu\rangle$. It follows that $\nabla_{U_i}A
= \langle \nabla_{U_i}A, U_i \rangle U_i + \langle \nabla_{U_i}A,
U_j \rangle U_j$, which allows to determine equation
(\ref{eqNablaUiA}) from equations (\ref{eqNablaUiUi}) and
(\ref{eqNablaUiUj}).

Since $U_i$ has constant length, we have $\langle\nabla_A
U_i,U_i\rangle=0$, and from Lemma \ref{thTwoEigenv} we get
$\langle\nabla_AU_i,A\rangle=-\langle\nabla_A A,U_i\rangle=0$ and
$\langle\nabla_AU_i,W_3\rangle=-\langle\nabla_A W_3,U_i\rangle=0$.
Using Lemma \ref{thThreeEigenv} and (\ref{eqNablaUiA}) we obtain
$0=\eR_{AU_iW_j\xi}
=(\lambda_i-\lambda_j)\langle\nabla_AU_i,W_j\rangle$ and hence
$\langle\nabla_AU_i,W_j\rangle=0$. Using this equality (with $i$
and $j$ interchanged) we get $0=A\langle W_i,J\xi\rangle
=\langle\enabla_AW_i,J\xi\rangle+\langle W_i,\enabla_AJ\xi\rangle
=b_i\langle\nabla_AW_i,U_i\rangle$, which yields
$\langle\nabla_AU_i,W_i\rangle = 0$. Thus we have $\nabla_AU_i =
\langle \nabla_AU_i,U_j \rangle U_j$. The latter inner product can
be calculated by using the explicit expression for $\bar{R}$,
Lemma \ref{thThreeEigenv} and (\ref{eqNablaUiA}) from
$$
\frac{(-1)^j}{4}(b_i^2-2b_j^2)  = \eR_{AU_i U_j \xi} =
(\lambda_i-\lambda_j)\langle\nabla_AU_i,U_j\rangle
-(\lambda_3-\lambda_j)(-1)^i\left(\lambda_i
+\frac{3b_i^2}{4(\lambda_3-\lambda_i)}\right).
$$
Altogether this now gives equation (\ref{eqNablaAUi}).

Since $A$ has constant length, we have
$\langle\nabla_AA,A\rangle=0$. Let $\nu \in \{1,2\}$. From
(\ref{eqNablaAUi}) we get $\langle\nabla_AA,U_\nu\rangle = 0$, and
from Lemma \ref{thTwoEigenv} we get $\langle\nabla_AA,W_\nu\rangle
= 0$. Next, we consider $0=A\langle JU_i,W_3\rangle
=\langle\enabla_AJU_i,W_3\rangle+\langle JU_i,\enabla_AW_3\rangle
= -\langle\nabla_AU_i,JW_3\rangle + \langle
JU_i,\enabla_AW_3\rangle$. The first term vanishes because of
(\ref{eqNablaAUi}), and in the second term we replace $JU_i$ by
its expression in terms of $A$ and $\xi$ to obtain $
0=\langle\nabla_AA,W_3\rangle$. This eventually implies equation
(\ref{eqNablaAA}).\end{proof}

\begin{corollary}\label{thAU1U2}
The integral curves of $A$ are geodesics in $M$ and the three
vector fields $A,U_1,U_2$ span an autoparallel distribution
${\mathfrak D}^\perp$, that is, ${\mathfrak D}^\perp$ is
integrable and its leaves are totally geodesic submanifolds of
$M$.
\end{corollary}

\begin{corollary}\label{thRelation2eigenvaluesbi}
The principal curvatures $\lambda_1,\lambda_2,\lambda_3$ and the
functions $b_1,b_2$ satisfy the equation
\[
0 = 3\big((\lambda_3-\lambda_2)^2b_1^2 +
(\lambda_3-\lambda_1)^2b_2^2\big) +
(\lambda_3-\lambda_1)(\lambda_3-\lambda_2)
\big(1+4\lambda_2(\lambda_3-\lambda_1)+
4\lambda_1(\lambda_3-\lambda_2)\big).
\]
\end{corollary}

\begin{proof} From Lemma \ref{thThreeEigenv} we get $1\! =
4\eR_{U_1U_2A\xi}
=4(\lambda_2-\lambda_3)\langle\nabla_{U_1}U_2,A\rangle
-4(\lambda_1-\lambda_3)\langle\nabla_{U_2}U_1,A\rangle$. The
assertion then follows by using equation
(\ref{eqNablaUiUj}).\end{proof}

\begin{lemma}\label{thRelation2eigenvalues}
If $i\in\{1,2\}$ and $m_i>1$, then $4\lambda_3\lambda_i=1$.
\end{lemma}

\begin{proof} Let $W_i\in\Gamma(T_{\lambda_i}\ominus\R U_i)$ be a
local unit vector field. Applying Lemma \ref{thEigenGauss} with
$X=W_i$ and $Y=A$, and taking into account (\ref{eqNablaAA}),
yields
\[
0 = 1-4\lambda_3\lambda_i
+8\langle\nabla_{W_i}A,\nabla_AW_i\rangle.
\]
We thus need to prove $\langle\nabla_{W_i}A,\nabla_AW_i\rangle =
0$. From (\ref{eqNablaAUi}) and (\ref{eqNablaAA}) we see that
$\nabla_AW_i \in \Gamma({\mathfrak D})$, and Lemma
\ref{thTwoEigenv} shows that $\nabla_AW_i$ is perpendicular to
$T_{\lambda_3} \cap {\mathfrak D}$. From Lemma \ref{thTwoEigenv}
we also see that $\nabla_{W_i}A$ is perpendicular to
$T_{\lambda_i} \cap {\mathfrak D}$. It thus suffices to prove that
$\langle\nabla_{W_i}A,W_j\rangle = 0$ for all $W_j \in
\Gamma(T_{\lambda_j} \cap {\mathfrak D})$, where $j \in \{1,2\}$
with $j \neq i$.

Let $\nu,\mu \in \{1,2\}$ with $\nu \neq \mu$. Then $0 =
W_i\langle U_\nu,JW_j\rangle =
\langle\enabla_{W_i}U_\nu,JW_j\rangle +\langle
U_\nu,\enabla_{W_i}JW_j\rangle =
\langle\nabla_{W_i}U_\nu,JW_j\rangle -\langle
JU_\nu,\enabla_{W_i}W_j\rangle$. As $JU_\nu=(-1)^\nu b_\mu A-b_\nu
\xi$, this implies
\begin{equation}\label{eqNablaWiWjAIntermediate}
\langle\nabla_{W_i}W_j,A\rangle=\frac{(-1)^\nu}{b_\mu}
\langle\nabla_{W_i}U_\nu,JW_j\rangle.
\end{equation}

As $T_{\lambda_j}$ is real, we can write
$JW_j=\widetilde{W}_i+\widetilde{W}_3$ with
$\widetilde{W}_i\in\Gamma(T_{\lambda_i}\ominus\R U_i)$ and
$\widetilde{W}_3\in\Gamma(T_{\lambda_3}\ominus\R A)$. We have
$\langle\nabla_{W_i}U_j,\widetilde{W}_i\rangle=0$ and
$0=W_i\langle\widetilde{W}_i,J\xi\rangle
=\langle\enabla_{W_i}\widetilde{W}_i,J\xi\rangle
+\langle\widetilde{W}_i,\enabla_{W_i}J\xi\rangle
=b_i\langle\nabla_{W_i}\widetilde{W}_i,U_i\rangle$ by Lemma
\ref{thTwoEigenv}, which implies
$\langle\nabla_{W_i}U_\nu,\widetilde{W}_i\rangle=0$ and hence
$\langle\nabla_{W_i}U_\nu,JW_j\rangle =
\langle\nabla_{W_i}U_\nu,\widetilde{W}_3\rangle$. From $0 =
W_i\langle\widetilde{W}_3,J\xi\rangle
=\langle\enabla_{W_i}\widetilde{W}_3,J\xi\rangle
+\langle\widetilde{W}_3,\enabla_{W_i}J\xi\rangle =
b_i\langle\enabla_{W_i}\widetilde{W}_3,U_i\rangle +
b_j\langle\enabla_{W_i}\widetilde{W}_3,U_j\rangle -
\lambda_i\langle JW_i,\widetilde{W}_j\rangle$ we obtain
$$
b_j\langle\nabla_{W_i}U_j,\widetilde{W}_3\rangle =
-\left(\frac{b_i^2}{2(\lambda_3-\lambda_i)}+\lambda_i\right)
\langle JW_i,\widetilde{W}_3\rangle
$$
by using Lemma \ref{thTwoEigenv}. From the same lemma it follows
that
\[
\langle\nabla_{W_i}U_i,\widetilde{W}_3\rangle
=\frac{b_i}{2(\lambda_3-\lambda_i)}\langle
JW_i,\widetilde{W}_3\rangle.
\]
Taking into account the last two equations,
(\ref{eqNablaWiWjAIntermediate}) becomes
\begin{equation}\label{eqNablaWiWjA}
\frac{(-1)^ib_i}{2b_j(\lambda_3-\lambda_i)}\langle
JW_i,\widetilde{W}_3\rangle = \langle\nabla_{W_i}W_j,A\rangle =
\frac{(-1)^i}{b_i\,b_j}\!\Bigl(\frac{b_i^2}{2(\lambda_3-\lambda_i)}
+\lambda_i\!\Bigr)\!\langle JW_i,\widetilde{W}_3\rangle.
\end{equation}
This readily implies $\lambda_i\langle
JW_i,\widetilde{W}_3\rangle=0$. Since at least one of the two
eigenvalues $\lambda_1,\lambda_2$ must be nonzero, it follows that
$\langle JW_1,\widetilde{W}_3\rangle = 0$ or $\langle
JW_1,\widetilde{W}_3\rangle = 0$. From (\ref{eqNablaWiWjA}) we
thus see that $\langle\nabla_{W_i}W_j,A\rangle = 0$ or
$\langle\nabla_{W_j}W_i,A\rangle = 0$. But from Lemma
\ref{thThreeEigenv} we know that
\[
0=\eR_{W_iW_jA\xi}
=(\lambda_j-\lambda_3)\langle\nabla_{W_i}W_j,A\rangle
-(\lambda_i-\lambda_3)\langle\nabla_{W_j}W_i,A\rangle,
\]
which implies that $\langle\nabla_{W_i}W_j,A\rangle = 0$ in both
cases. This finishes the proof.
\end{proof}

From Lemma \ref{thRelation2eigenvalues} we immediately get

\begin{corollary}\label{multonelambda}
$m_1=1$ or $m_2=1$.
\end{corollary}

According to Corollary \ref{multonelambda} we may assume that
$m_2=1$, that is, $T_{\lambda_2} = \R U_2$. We will now
distinguish the two cases $m_1>1$ and $m_1=1$.\medskip

{\it Case 1: $m_1>1$.} Then we have $4\lambda_1\lambda_3=1$ by
Lemma \ref{thRelation2eigenvalues} and $J(T_{\lambda_1}\ominus\R
U_1)\subset T_{\lambda_3}\ominus\R A$ because $T_{\lambda_1}$ is
real. Let $W_1\in\Gamma(T_{\lambda_1}\ominus\R U_1)$ and
$W_3\in\Gamma(T_{\lambda_3}\ominus\R A)$. Then $0 =W_1\langle
W_3,J\xi\rangle=\langle\enabla_{W_1}W_3,J\xi\rangle +\langle
W_3,\enabla_{W_1}J\xi\rangle =
b_1\langle\nabla_{W_1}W_3,U_1\rangle +
b_2\langle\nabla_{W_1}W_3,U_2\rangle - \lambda_1\langle
JW_1,W_3\rangle$. Using Lemma \ref{thTwoEigenv} this implies
\begin{equation}\label{thNablaW1U2W3}
\langle\nabla_{W_1}U_2,W_3\rangle =
-\frac{1}{b_2}\left(\frac{b_1^2}{2(\lambda_3-\lambda_1)}+\lambda_1\right)
\langle JW_1,W_3\rangle.
\end{equation}
Next, using $JU_2 = b_1A-b_2\xi$, we get $0=W_1\langle
W_1,JU_2\rangle =\langle\enabla_{W_1}W_1,JU_2\rangle +\langle
W_1,\enabla_{W_1}JU_2\rangle = \langle\nabla_{W_1}U_2,JW_1\rangle
+b_1\langle W_1,\nabla_{W_1}A\rangle - b_2\langle
W_1,\enabla_{W_1}\xi\rangle$. We now assume that $W_1$ has length
one. Using Corollary \ref{thNablaXX} this implies
\begin{equation}\label{thNablaW1U2JW1}
\langle\nabla_{W_1}U_2,JW_1\rangle = -b_2\lambda_1.
\end{equation}
Comparing (\ref{thNablaW1U2W3}) with $W_3 = JW_1$ and
(\ref{thNablaW1U2JW1}), and using $b_1^2+b_2^2=1$, implies
$2\lambda_1(\lambda_1-\lambda_3)=1$. Together with
$4\lambda_1\lambda_3=1$ this implies $\lambda_1=\sqrt{3}/2$ and
$\lambda_3=\sqrt{3}/6$, where we assume that the orientation of
$\xi$ is such that $\lambda_3 > 0$.

We now apply Lemma \ref{thEigenGauss} with $X=W_1$ and $Y=U_2$,
and use Corollary \ref{thNablaXX} and (\ref{eqNablaUiUi}), to
obtain
\begin{equation}\label{tharteq}
0 = (\lambda_2-\lambda_1)(1-4\lambda_1\lambda_2
+8\langle\nabla_{W_1}U_2,\nabla_{U_2}W_1\rangle) -b_2\langle
\nabla_{W_1}U_2,JW_1\rangle +2b_2\langle
\nabla_{U_2}W_1,JW_1\rangle.
\end{equation}
Using Lemma \ref{thTwoEigenv} we easily get
$\nabla_{W_1}U_2\in\Gamma(T_{\lambda_3})$, and (\ref{eqNablaUiA})
shows that $\langle \nabla_{U_2}W_1,A\rangle = 0$. From
(\ref{thNablaW1U2W3}) we thus get
\begin{equation}\label{thNablaW1U2NablaU2W1}
\langle\nabla_{W_1}U_2,\nabla_{U_2}W_1\rangle
=\langle\nabla_{W_1}U_2,JW_1\rangle\langle\nabla_{U_2}W_1,JW_1\rangle.
\end{equation}
From Lemma \ref{thThreeEigenv} and (\ref{thNablaW1U2JW1}) we
obtain $b_2=4\eR_{U_2W_1JW_1\xi}
=4(\lambda_1-\lambda_3)\langle\nabla_{U_2}W_1,JW_1\rangle
+4(\lambda_2-\lambda_3)b_2\lambda_1$ and hence
\begin{equation}\label{thNablaU2W1JW1}
\langle\nabla_{U_2}W_1,JW_1\rangle =
\frac{b_2}{4(\lambda_1-\lambda_3)}(1-4\lambda_1(\lambda_2-\lambda_3))
= b_2\lambda_1(1-2\lambda_1\lambda_2).
\end{equation}
Inserting (\ref{thNablaW1U2W3}), (\ref{thNablaW1U2NablaU2W1}) and
(\ref{thNablaU2W1JW1}) into (\ref{tharteq}) yields
$$
0 = 12(3b_2^2-1)\lambda_2^2 + 4\sqrt{3}(2-9b_2^2)\lambda_2 +
3(9b_2^2-1).
$$
On the other hand, inserting the above particular values for
$\lambda_1$ and $\lambda_3$ into the equation in Corollary
\ref{thRelation2eigenvaluesbi}, and replacing $b_1^2$ by
$1-b_2^2$, yields
$$
0 = 12(9b_2^2+1)\lambda_2^2 - 4\sqrt{3}(2+9b_2^2)\lambda_2 -
3(9b_2^2-1).
$$
Adding up the previous two equations gives
$\lambda_2(2\lambda_2-\sqrt{3})=0$. As $\lambda_2 \neq \lambda_1 =
\sqrt{3}/2$ we therefore get $\lambda_2 = 0$, which implies $b_2^2
= 1/9$ and $b_1^2 = 8/9$.

\medskip
{\it Case 2: $m_1=1$.} In this case we have $T_{\lambda_1}=\R
U_1$, and   $T_{\lambda_3}\ominus\R A = {\mathfrak D}$ is a
$J$-invariant distribution. Let $W_3 \in \Gamma({\mathfrak D})$ be
of unit length. We have $0 = W_3\langle JW_3,J\xi\rangle
=\langle\enabla_{W_3}JW_3,J\xi\rangle +\langle
JW_3,\enabla_{W_3}J\xi\rangle$, and applying Lemma
\ref{thTwoEigenv} this yields
\begin{equation}\label{equationforb1b2}
(\lambda_3-\lambda_2)b_1^2 + (\lambda_3-\lambda_1)b_2^2+
4\lambda_3(\lambda_3-\lambda_1)(\lambda_3-\lambda_2) = 0.
\end{equation}
Together with $b_1^2+b_2^2=1$ this implies
\begin{equation}\label{eqbi2}
b_i^2=\frac{\lambda_3-\lambda_i}{\lambda_j-\lambda_i}
(1+4\lambda_3(\lambda_3-\lambda_j))\qquad (i,j \in \{1,2\},\ i
\neq j).
\end{equation}
Inserting these expressions for $b_1^2$ and $b_2^2$ into the
equation in Corollary \ref{thRelation2eigenvaluesbi} yields
\begin{equation}\label{equation1}
(\lambda_1-\lambda_2)^2 - (\lambda_1+\lambda_2-4\lambda_3)^2 =
1-4\lambda_3^2.
\end{equation}
We now apply Lemma \ref{thEigenGauss} with $X = W_3$ and $Y=U_i$,
which gives
\begin{equation}\label{GausseqnW3Ui}
\begin{array}{rcl}
0&=&({\lambda_i}-{\lambda_3})\big( 1-4{\lambda_i}\,{\lambda_3}
                +8\langle\nabla_{W_3}U_i,\nabla_{U_i} W_3\rangle
        -4\langle\nabla_{W_3}W_3,\nabla_{U_i}U_i\rangle\big)\\
&&   -b_i\langle \nabla_{W_3} U_i,JW_3\rangle + 2b_i\langle
\nabla_{U_i} W_3,JW_3\rangle.
\end{array}
\end{equation}
Let $i,j \in \{1,2\}$ and $i \neq j$. Then we have $0 = W_3\langle
U_i,J\xi\rangle =\langle\enabla_{W_3}U_i,J\xi\rangle +\langle
U_i,\enabla_{W_3}J\xi\rangle
=b_j\langle\nabla_{W_3}U_i,U_j\rangle$ and hence
$\langle\nabla_{W_3}U_i,U_j\rangle = 0$. For $\widetilde{W}_3 \in
\Gamma({\mathfrak D})$ we have $4(\lambda_3-\lambda_i)\langle
\nabla_{W_3}U_i,\widetilde{W}_3\rangle =
4(\lambda_i-\lambda_3)\langle
\nabla_{W_3}\widetilde{W}_3,U_i\rangle = b_i\langle
JW_3,\widetilde{W}_3\rangle$ and $\langle \nabla_{W_3}U_i,A\rangle
= -\langle \nabla_{W_3}A,U_i\rangle = 0$ by Lemma
\ref{thTwoEigenv}. Altogether this gives
$4(\lambda_3-\lambda_i)\nabla_{W_3}U_i = b_iJW_3$, and together
with (\ref{eqNablaUiUi}) equation (\ref{GausseqnW3Ui}) now becomes
$$
0 = 4(\lambda_3-\lambda_i)^2(1-4\lambda_i\lambda_3) -
12(-1)^ib_1b_2(\lambda_3-\lambda_i)\langle \nabla_{W_3}W_3 , A
\rangle + b_i^2.
$$
Multiplying this equation with $\lambda_3-\lambda_j$, then adding
the two equations for $i=1$ and $i=2$, and then using
(\ref{equationforb1b2}) yields
\begin{equation}\label{equation2}
4\lambda_3(1+\lambda_1^2+\lambda_2^2) -
(\lambda_1+\lambda_2)(1+4\lambda_3^2) = 0.
\end{equation}
If $\lambda_3=0$, we immediately get $\lambda_1,\lambda_2 \in
\{\pm 1/2\}$ from (\ref{equation1}) and (\ref{equation2}). From
now on we assume $\lambda_3 \neq 0$. If we put
$x=\lambda_1-\lambda_2$ and $y = \lambda_1+\lambda_2-4\lambda_3$,
equations (\ref{equation1}) and (\ref{equation2}) are equivalent
to
$$
x^2 - y^2 = 1-4\lambda_3^2\ ,\ x^2 +
\left(y-\frac{1-12\lambda_3^2}{4\lambda_3}\right)^2 =
\frac{1+16\lambda_3^4}{16\lambda_3^2}.
$$
Obviously, these are the equations of a hyperbola and a circle. It
is straightforward to calculate their common points, namely
$$
(x,y) = \left(\pm \sqrt{1-3\lambda_3^2}\ ,\ -\lambda_3\right)\ \
{\rm and}\ \ (x,y) = \left(\pm\frac{1}{4\lambda_3}\ ,\
\frac{1-8\lambda_3^2}{4\lambda_3}\right)\ ,
$$
where the first possibility only arises if $3\lambda_3^2 \leq 1$.
Taking into account that $\lambda_1$ and $\lambda_2$ are different
from $\lambda_3$, this eventually implies
\begin{equation}\label{thlambda12}
\lambda_1 = \frac{1}{2}\left(3\lambda_3 -
\sqrt{1-3\lambda_3^2}\right)\ ,\ \lambda_2 =
\frac{1}{2}\left(3\lambda_3 + \sqrt{1-3\lambda_3^2}\right)
\end{equation}
where we assume without loss of generality that $\lambda_1 <
\lambda_2$. Obviously, we get a solution only if $3\lambda_3^2
\leq 1$. If $|\lambda_3|=1/2$ or $|\lambda_3|=1/\sqrt{3}$, then
the three principal curvatures cannot be different. Suppose  that
$1/2<|\lambda_3|<1/\sqrt{3}$. From (\ref{equationforb1b2}) and
(\ref{thlambda12}) we get
\[
\frac{b_1^2}{2\lambda_3(\lambda_3-\sqrt{1-3\lambda_3^2})}
+\frac{b_2^2}{2\lambda_3(\lambda_3+\sqrt{1-3\lambda_3^2})}=1.
\]
If $1/2<|\lambda_3|<1/\sqrt{3}$, elementary calculations
show that $0<2\lambda_3(\lambda_3-\sqrt{1-3\lambda_3^2})<1$ and
$0<2\lambda_3(\lambda_3+\sqrt{1-3\lambda_3^2})<1$. Therefore the
last equation is the equation of an ellipse centered at the origin
and with axes of length less than 1. Obviously such an ellipse has
no points of intersection with the circle $b_1^2+b_2^2=1$. This
shows that $|\lambda_3|<1/2$.

We summarize the discussion in this section in

\begin{theorem}\label{thStep1}
Let $M$ be a connected real hypersurface in ${\mathbb C}H^n$, $n
\geq 3$, with three distinct constant principal curvatures
$\lambda_1,\lambda_2,\lambda_3$, and suppose that $M$ is not a
Hopf hypersurface. Then, with a suitable labelling of the
principal curvatures, we have $J\xi = b_1U_1 + b_2U_2$ with some
real numbers $b_1,b_2 >0$, where $U_i$ denotes the orthogonal
projection of $J\xi$ onto $T_{\lambda_i}$ normalized to length
one. There exists a unit vector field $A \in
\Gamma(T_{\lambda_3})$ such that $JA = b_2U_1 - b_1U_2$. The
subbundle $\R U_1 \oplus \R U_2$ is real, and the subbundle $\R A
\oplus \R U_1 \oplus \R U_2 \oplus \R\xi$ is complex. Moreover,
$m_2 = 1$ and one the following two cases holds:
\begin{enumerate}
\item[(i)] $m_1 > 1$, $\lambda_1=\sqrt{3}/2$, $\lambda_2=0$,
$\lambda_3=\sqrt{3}/6$, $b_1=2\sqrt{2}/3$, $b_2=1/3$, the
subbundle $T_{\lambda_1}\ominus\R U_1$ is real, and
$J(T_{\lambda_1}\ominus\R U_1)\subset T_{\lambda_3}$. \item[(ii)]
$m_1 = 1$, $-1/2 < \lambda_3 < 1/2$,
$\lambda_1=\frac{1}{2}(3\lambda_3-\sqrt{1-3\lambda_3^2})$,
$\lambda_2=\frac{1}{2}(3\lambda_3+\sqrt{1-3\lambda_3^2})$, and
$$
b_i^2=\frac{\lambda_3-\lambda_i}{\lambda_j-\lambda_i}
(1+4\lambda_3(\lambda_3-\lambda_j))\qquad (i,j \in \{1,2\},\ i
\neq j).
$$
\end{enumerate}
\end{theorem}

\section{Proof of Theorem \ref{thClassification}}\label{proofmain}

In this section we prove Theorem \ref{thClassification}. Let $M$
be a connected real hypersurface in $\CH$ with three distinct
constant principal curvatures. If $M$ is a Hopf hypersurface, it
was shown in \cite{Be89} that $M$ is an open part of a tube around
a totally geodesic $\CH[k]\subset\CH$ for some
$k\in\{1,\dots,n-2\}$, or of a tube with radius $r \neq
\ln(2+\sqrt{3})$ around a totally geodesic $\R H^n\subset\CH$. We
can therefore assume that $M$ is not a Hopf hypersurface. Then $M$
must satisfy one of the two possibilities described in Theorem
\ref{thStep1}. The result will follow from a thorough
investigation of the possible focal sets and equidistant
hypersurfaces of $M$ by means of Jacobi field theory.

For $r\in\R$ we define the smooth map $\Phi^r : M\to \CH,\ p
\mapsto \Phi^r(p)=\exp_p(r\xi_p)$, where $\exp_p$ is the
exponential map of $\CH$ at $p$. Geometrically this means that we
assign to $p$ the point in $\CH$ which is obtained by travelling
for the distance $r$ along the geodesic $c_p(t)=\exp_p(t\xi_p)$ in
direction of the normal vector $\xi_p$ (for $r >0$; for $r<0$ one
sets off in direction $-\xi_p$; and for $r=0$ there is no movement
at all). For $v \in T_pM$ we denote by $B_v$ the parallel vector
field along the geodesic $c_p$ with $B_v(0) = v$, and by $\zeta_v$
the Jacobi field along $c_p$ with $\zeta_v(0)=v$ and
$\zeta_v'(0)=-S_pv$. Note that $\zeta_v$ is the unique solution of
the linear differential equation
$$
4\zeta_v''-\zeta_v-3\langle \zeta_v,J\dot{c}_p\rangle
J\dot{c}_p=0\ ,\ \zeta_v(0)=v\ ,\ \zeta_v'(0)=-S_pv,
$$
where $\dot{c}_p$ denotes the tangent vector field of $c_p$ and
the prime $'$ indicates the covariant derivative of a vector field
along $c_p$. For $v \in T_{\lambda_i}(p)$ we have the explicit
expression
$$
\zeta_v(t)=f_i(t)B_v(t)+\langle v,J\xi\rangle g_i(t)J\dot{c}_p(t)
$$
with
$$
\begin{array}{rcl}
f_i(t)    &=& \displaystyle
        \cosh(t/2)-2\lambda_i\sinh(t/2),\\
\noalign{\medskip} g_i(t)    &=&
\displaystyle\left(\cosh(t/2)-1\right)
    \left(1+2\cosh(t/2)-2\lambda_i\sinh(t/2)\right).
\end{array}
$$
Finally, we define a vector field $\eta^r$ along the map $\Phi^r$
by $\eta^r_p=\dot{c}_p(r)$. The relation between the map $\Phi^r$,
the vector field $\eta^r$ and the Jacobi field $\zeta_v$ is given
by
$$\zeta_v(r) = \Phi^r_*v\ ,\
\zeta_v'(r) = \enabla_{v}\eta^r,
$$
where $\Phi^r_*$ denotes the differential of $\Phi^r$. The
singularities of $\Phi^r$ are focal points of $M$ and can be
calculated using Jacobi fields from the equation $\zeta_v(r) =
\Phi^r_*v$. We will see that in case (i) of Theorem \ref{thStep1}
there exists a particular distance $r$ at which the map $\Phi^r$
has constant rank $2n-m_1$, which means that the image of $\Phi^r$
forms locally a submanifold of codimension $m_1$. In case (ii) of
Theorem \ref{thStep1} there exists a particular distance $r$ at
which the map $\Phi^r$ has constant rank $2n-1$ and the image is
locally a minimal real hypersurface. We then use the equation
$\zeta_v'(r) = \enabla_{v}\eta^r$ to obtain some information about
the second fundamental form of these submanifolds. We continue
using the notation introduced in Section \ref{princurv}.

\medskip
{\it Case 1: $m_1>1$.} We define $u_i = (U_i)_p$ and $r =
\ln(2+\sqrt{3})$. For $v \in T_pM$ we denote by $v_i$ the
orthogonal projection of $v$ onto $T_{\lambda_i}(p)$. Using the
equation $\Phi^r_*v = \zeta_v(r)$ and the explicit expression for
the Jacobi fields, we obtain
\begin{eqnarray*}
9\Phi^r_*v &=& 3\sqrt{6}B_{v_3}(r)
    +\big(4\langle v_1,u_1\rangle
    +(4\sqrt{2}-2\sqrt{3})\langle v_2,u_2\rangle\big)B_{u_1}(r)\\
&&  +\big(\sqrt{2}\langle v_1,u_1\rangle
        +(2+4\sqrt{6})\langle v_2,u_2\rangle\big)B_{u_2}(r).
\end{eqnarray*}
This shows that $\Phi^r_*v=0$ if and only if $v \in
T_{\lambda_1}(p)\ominus\R u_1$.  Therefore the rank of $\Phi^r$ is
constant and equal to $2n-m_1$. This means that for every point in
$M$ there exists an open neighborhood $\mathcal{V}$ such that
$\mathcal{W}=\Phi^r(\mathcal{V})$ is an embedded submanifold of
$\CH$ and $\Phi^r:\mathcal{V}\to\mathcal{W}$ is a submersion. Let
$p \in {\mathcal V}$ and $q = \Phi^r(p) \in {\mathcal W}$. The
above expression for the differential of $\Phi^r$ shows that the
tangent space $T_q{\mathcal W}$ of ${\mathcal W}$ at $q$ is
obtained by parallel translation of $T_{\lambda_3}(p) \oplus \R
u_1 \oplus \R u_2$ along the geodesic $c_p$ from $p = c_p(0)$ to
$q = c_p(r)$. Hence, the normal space $\nu_q{\mathcal W}$ of
${\mathcal W}$ at $q$ is obtained by parallel translation of $\R
\xi_p \oplus (T_{\lambda_1}(p) \ominus \R u_1)$ along $c_p$ from
$p$ to $q$. This shows in particular that ${\mathcal W}$ has
totally real normal bundle.

Clearly, $\eta^r_p = B_{\xi_p}(r)$ is a unit normal vector of
${\mathcal W}$ at $q$. For the shape operator $S^r$ of ${\mathcal
W}$ we have $S^r_{\eta^r_p}\Phi^r_*v = -(\enabla_{v}\eta^r)^\top =
-(\zeta_v'(r))^\top$, where $(\cdot )^\top$ denotes the component
tangent to ${\mathcal W}$. Using the explicit expression for the
Jacobi fields we easily get
\begin{equation}\label{thSretarpBv3}
S^r_{\eta^r_p}B_{v_3}(r) = 0 \ {\rm for\ all}\ v_3 \in
T_{\lambda_3}(p).
\end{equation}
Moreover, $S^r_{\eta^r_p}$ leaves
$\R B_{u_1}(r) \oplus \R B_{u_2}(r)$ invariant and has the matrix
representation
$$
\frac{1}{18} \left(\begin{array}{cc}
4\sqrt{2} &   -7\\
-7     &   -4\sqrt{2}
\end{array}\right),
$$
with respect to $B_{u_1}(r),B_{u_2}(r)$. Since $3JA_p = u_1 -
2\sqrt{2}u_2$ and $3J\xi_p = 2\sqrt{2}u_1 + u_2$, the above matrix
representation yields
\begin{equation}\label{thSretarpBJa}
2S^r_{\eta^r_p}B_{JA_p}(r) = J\eta^r_p\ ,\
2S^r_{\eta^r_p}J\eta^r_p = B_{JA_p}(r). \end{equation} As
$J(\nu_q{\mathcal W} \ominus \R \eta^r_p)$ is contained in the
parallel translate of $T_{\lambda_3}(p)$ along $c_p$ from $p$ to
$q$, (\ref{thSretarpBv3}) and the linearity of $S^r_{\eta^r_p}$
show that
\begin{equation}\label{eqSretarp}
2S^r_{\eta^r_p}J\tilde{\eta} = \langle\eta^r_p,\tilde{\eta}\rangle
B_{JA_p}(r)\ {\rm for\ all}\ \tilde{\eta} \in \nu_q{\mathcal W}.
\end{equation}
As a special case we get $S^r_{\eta^r_p}J\tilde{\eta} = 0$ for all
$p \in {\mathcal V}$ and $\tilde{\eta} \in \nu_q{\mathcal W}
\ominus \R \eta^r_p$. From the Gau\ss\ formula and $\enabla J = 0$
one easily gets $S^r_{\tilde{\eta}}J\eta^r_p =
S^r_{\eta^r_p}J\tilde{\eta}$ and hence
\begin{equation}\label{thSretaJeta}
S^r_{\tilde{\eta}}J\eta^r_p = 0\ {\rm for\ all}\ \tilde{\eta} \in
\nu_q{\mathcal W} \ominus \R \eta^r_p.
\end{equation}
Now let $\gamma$ be any curve in $(\Phi^r)^{-1}(\{q\}) \cap
{\mathcal V}$ with $\gamma(0) = p$. Since $\eta^r_p$ and
$\eta^r_{\gamma(t)} - \langle \eta^r_{\gamma(t)} , \eta^r_p
\rangle \eta^r_p$ are perpendicular, (\ref{thSretaJeta}), the
linearity of $\eta \mapsto S^r_\eta$ and (\ref{eqSretarp}) imply
$$
0 = 2S_{\eta^r_{\gamma(t)} - \langle \eta^r_{\gamma(t)} , \eta^r_p
\rangle \eta^r_p}J\eta^r_p = 2S_{\eta^r_{\gamma(t)}}J\eta^r_p -
\langle \eta^r_{\gamma(t)} , \eta^r_p \rangle B_{JA_p}(r).
$$
On the other hand, (\ref{eqSretarp}) with $\gamma(t)$ instead of
$p$ gives
$$2S^r_{\eta^r_{\gamma(t)}}J\eta^r_p =
\langle\eta^r_{\gamma(t)},\eta^r_p\rangle B_{JA_{\gamma(t)}}(r).$$
The previous two equations show that the map $\tilde{p} \mapsto
B_{JA_{\tilde{p}}}(r)$ is of constant value $z \in T_q{\mathcal
W}$ on the connected component ${\mathcal V}_o$ of
$(\Phi^r)^{-1}(\{q\}) \cap {\mathcal V}$ containing $p$. Note that
$z$ has length one because of $z = B_{JA_p}(r)$. For all $v_1 \in
T_{\lambda_1}(p) \ominus \R u_1$ we have $\enabla_{v_1}\eta^r =
\zeta_{v_1}'(r) = (-1/\sqrt{2})B_{v_1}(r)$, which implies that
$\eta^r$ is a local diffeomorphism from ${\mathcal V}_o$ into the
unit sphere in $\nu_q{\mathcal W}$. Thus $\eta^r({\mathcal V}_o)$
is an open subset of the unit sphere in $\nu_q{\mathcal W}$. Since
$S^r_\eta$ depends analytically on $\eta \in \nu_q{\mathcal W}$,
we conclude from (\ref{thSretarpBv3}) and (\ref{thSretarpBJa})
that
$$
2S^r_\eta J\eta = z\ ,\ 2S^r_\eta z = J\eta\ ,\ S^r_\eta v = 0\
{\rm for\ all}\ \eta \in \nu_q{\mathcal W}, v \in T_q{\mathcal W}
\ominus J(\nu_q{\mathcal W} \ominus \R\eta) \ominus \R z.
$$
Therefore the second fundamental form $I\!I^r_q$ of ${\mathcal W}$
at $q$ is given by the trivial bilinear extension of
$2I\!I^r_q(z,J\eta) = \eta$ for all $\eta \in \nu_q{\mathcal W}$.
The construction of $z$ shows that it depends smoothly on the
point $q \in {\mathcal W}$. Hence there exists a unit vector field
$Z$ on ${\mathcal W}$ such that the second fundamental form
$I\!I^r$ of ${\mathcal W}$ is given by the trivial bilinear
extension of $2I\!I^r(Z,J\eta) = \eta$ for all $\eta \in
\Gamma(\nu{\mathcal W})$. From Theorem \ref{ruledk} we see that
${\mathcal W}$ is holomorphically congruent to an open part of the
ruled minimal submanifold $W^{2n-m_1}$. Thus we have proved that
locally $M$ lies on a tube with radius $r = \ln(2+\sqrt{3})$
around a ruled minimal submanifold holomorphically congruent to
$W^{2n-m_1}$. This finally implies that $M$ is holomorphically
congruent to an open part of the tube with radius $r =
\ln(2+\sqrt{3})$ around $W^{2n-m_1}$.

\medskip
{\it Case 2: $m_1=1$.} If $\lambda_3= 0$, then $\lambda_1 = -1/2$
and $\lambda_2 = 1/2$, and it follows from Theorem \ref{ruled1}
that $M$ is holomorphically congruent to an open part of the ruled
minimal hypersurface $W^{2n-1}$. If $0 < |\lambda_3| < 1/2$, we
can write $2\lambda_3 = \tanh(r/2)$ with some $0 \neq r \in \R$.

Let $p \in M$ and define $u_i = (U_i)_p$. Using the equation
$\Phi^r_*v = \zeta_v(r)$ and the explicit expression for the
Jacobi fields, we obtain
$$
\Phi^r_* v_3 = \sech(r/2)B_{v_3}(r)\ {\rm for\ all}\ v_3 \in
T_{\lambda_3}(p)
$$
and
$$
\left( \begin{array}{@{}c@{}} \Phi^r_*u_1 \\ \Phi^r_*u_2
\end{array} \right) = D(r) \left( \begin{array}{@{}c@{}} B_{u_1}(r) \\
B_{u_2}(r)
\end{array} \right)
$$
with
$$
D(t) = \left(
\begin{array}{@{}cc@{}}
f_1(t)+b_1^2 g_1(t)   &   b_1b_2g_1(t)\\
b_1b_2g_2(t)  & f_2(t)+b_2^2 g_2(t)
\end{array}\right).
$$
As $\det(D(r)) = \sech^3(r/2)$, we can now conclude that
$\Phi^r_*$ has maximal rank everywhere. This means that for every
point in $M$ there exists an open neighborhood $\mathcal{V}$ such
that $\mathcal{W}=\Phi^r(\mathcal{V})$ is an embedded real
hypersurface of $\CH$ and $\Phi^r:\mathcal{V}\to\mathcal{W}$ is a
diffeomorphism. Let $p \in {\mathcal V}$ and $q = \Phi^r(p) \in
{\mathcal W}$. The tangent space $T_q{\mathcal W}$ of ${\mathcal
W}$ at $q$ is obtained by parallel translation of $T_p{\mathcal
V}$ along the geodesic $c_p$ from $p = c_p(0)$ to $q = c_p(r)$,
and $\eta^r_p$ is a unit normal vector of ${\mathcal W}$ at $q$.

For the shape operator $S^r$ of ${\mathcal W}$ we have
$S^r_{\eta^r_p}\Phi^r_*v = -\enabla_{v}\eta^r = -\zeta_v'(r)$.
Since $f_3'(r) = 0$ we immediately get
$$
S^r_{\eta^r_p}B_{v_3}(r) = 0 \ \mbox{for all}\ v_3 \in
T_{\lambda_3}(p),
$$
and for $\Phi^r_*u_1$ and $\Phi^r_*u_2$ we get
$$
\left( \begin{array}{@{}c@{}} S^r_{\eta^r_p}\Phi^r_*u_1 \\
S^r_{\eta^r_p}\Phi^r_*u_2 \end{array}\right)
= C(r) \left( \begin{array}{@{}c@{}} B_{u_1}(r) \\
B_{u_2}(r)
\end{array} \right)
$$
with $C(r) = -D'(r)D(r)^{-1}$. A tedious calculation shows that
$\det(D'(r)) = -\sech^3(r/2)/4$ and $(\det(D))'(r) = 0$, which
implies
$$\det(C(r)) = \frac{\det(D'(r))}{\det(D(r))} = -\frac{1}{4}
\quad \mbox{and}\quad  \tr(C(r)) =
-\frac{(\det(D))'(r)}{\det(D(r))} = 0.$$ From this we easily see
that the eigenvalues of $C(r)$ are $\pm 1/2$. Altogether we now
get that ${\mathcal W}$ has three distinct constant principal
curvatures $0$, $+1/2$ and $-1/2$ with corresponding
multiplicities $2n-3$, $1$ and $1$, respectively. It follows from
Theorem \ref{ruled1} that ${\mathcal W}$ is holomorphically
congruent to an open part of the ruled real hypersurface
$W^{2n-1}$. From this we eventually conclude that $M$ is
holomorphically congruent to an open part of an equidistant
hypersurface to $W^{2n-1}$.

This finishes the proof of Theorem \ref{thClassification}.


\bigskip
\noindent {\sc Department of Mathematics, University College,
Cork, Ireland}\\
\noindent {\sc Email:} j.berndt@ucc.ie

\smallskip
\noindent {\sc Department of Geometry and Topology, Faculty of
Mathematics,\\ University of Santiago de Compostela, Spain}\\
\noindent {\sc Email:} xtjosec@usc.es


\begin{thebibliography}{9999}

\bibitem{Be89}
J.\ Berndt: Real hypersurfaces with constant principal curvatures
in complex hyperbolic space, {\it J.\ Reine Angew.\ Math.} {\bf
395} (1989), 132--141.

\bibitem{Be98}
J.\ Berndt: {Homogeneous hypersurfaces in hyperbolic spaces}, {\it
Math.\ Z.} {\bf 229} (1998), 589--600.

\bibitem{BB01}
J.\ Berndt, M.\ Br{\"u}ck: Cohomogeneity one actions on hyperbolic
spaces, {\it J.\ Reine Angew.\ Math.} {\bf 541} (2001), 209--235.

\bibitem{BCO03}
J.\ Berndt, S.\ Console, C.\ Olmos: {\it Submanifolds and
holonomy}, Chapman \& Hall/CRC Research Notes in Mathematics {\bf
434}, Chapman \& Hall/CRC, Boca Raton, FL, 2003.

\bibitem{BT04}
J.\ Berndt, H.\ Tamaru: Cohomogeneity one actions on noncompact
symmetric spaces of rank one, arXiv:math/0505490.

\bibitem{BTV95}
J.\ Berndt, F.\ Tricerri, L.\ Vanhecke: {\it Generalized
Heisenberg groups and Damek--Ricci harmonic spaces}, Lecture Notes
in Mathematics {\bf 1598}, Springer--Verlag, Berlin, 1995.

\bibitem{Ca38}
\'{E}.\ Cartan: Familles de surfaces isoparam\'etriques dans les
espaces \`{a} courbure constante, {\it Ann.\ Mat.\ Pura Appl., IV.
Ser.} {\bf 17} (1938), 177--191.

\bibitem{Mo85}
S.\ Montiel: Real hypersurfaces of a complex hyperbolic space,
{\it J.\ Math.\ Soc.\ Japan} {\bf 37} (1985), 515--535.

\bibitem{Sa99}
J.\ Saito: Real hypersurfaces in a complex hyperbolic space with
three constant principal curvatures, {\it Tsukuba J.\ Math.} {\bf
23} (1999), 353--367.

\bibitem{Th00} G.\ Thorbergsson: A survey on isoparametric
hypersurfaces and their generalizations, {\it Handbook of
differential geometry, Vol. I}, 963--995, North-Holland,
Amsterdam, 2000.


\end{thebibliography}
\end{document}